\documentclass{amsart}

\usepackage{amssymb}

\newtheorem{thm}{Theorem}
\newtheorem{lem}[thm]{Lemma}
\newtheorem{prop}[thm]{Proposition}
\newtheorem{cor}[thm]{Corollary}

\theoremstyle{definition}
\newtheorem{defn}[thm]{Definition}
\newtheorem{exmp}[thm]{Example}
\newtheorem{rem}[thm]{Remark}

\numberwithin{equation}{subsection}
\numberwithin{thm}{subsection}

\newcommand\idd{\mathop {\fam 0 id}\nolimits}
\newcommand\Id{{\fam 0 I}}
\newcommand\Curr{\mathop {\fam 0 Cur}\nolimits}
\newcommand\Ann{\mathop {\fam 0 Ann}\nolimits}

\newcommand\Cend{\mathop {\fam 0 Cend}\nolimits}
\newcommand\gc{\mathop {\fam 0 gc}\nolimits}
\newcommand\oc{\mathop {\fam 0 oc}\nolimits}
\newcommand\spc{\mathop {\fam 0 spc}\nolimits}

\newcommand\Gll{\mathop {\fam 0 gl}\nolimits}
\newcommand\Sll{\mathop {\fam 0 sl}\nolimits}

\newcommand\Kerr{\mathop {\fam 0 Ker}\nolimits}
\newcommand\diag{\mathop {\fam 0 diag}\nolimits}
\newcommand\ideal{\mathrel {\triangleleft}}
\newcommand\Coeff{\mathop {\fam 0 Coeff}}
\newcommand\Diff{\mathop {\fam 0 Diff}\nolimits}

\newcommand{\oo}[1]{\mathrel{{\circ }_{#1}} }

\newcommand\Ress{\mathop{ \fam 0 Res}}
\newcommand\GKdim{\mathop {\fam 0 GKdim}\nolimits}

\newcommand\End{\mathop {\fam 0 End}\nolimits}
\newcommand\Zset{{\mathbb Z}}

\begin{document}

\title{Associative conformal algebras
 with finite faithful representation}

\thanks{Partially supported by RFBR  and SSc-2269.2003}

\author{Pavel Kolesnikov}
\email{pavelsk@kias.re.kr}

\address{Korea Institute for Advanced Study}

\begin{abstract}
We describe irreducible conformal subalgebras of $\Cend_N$
and build the structure theory of associative conformal algebras with
finite faithful representation.
\end{abstract}



\maketitle


\section{Introduction}\label{sec1}

The notion of a conformal algebra appears as an algebraic
language describing
singular part of operator product expansion (OPE) in conformal
field theory~\cite{3}.
The explicit algebraic exposition of this
theory
(see, e.g.,
\cite{7,14,18})
  leads to the notion of a vertex (chiral) algebra.
Roughly speaking, conformal algebras correspond to vertex algebras
 by the same way as Lie algebras correspond to their associative
enveloping algebras (see~\cite{BK,20} for a detailed explanation).

In the recent years a great advance in the structure theory
 of associative and Lie conformal algebras
 of finite type has been obtained.
In \cite{11,12,16,17,21},
simple and semisimple Lie conformal (super)algebras of finite type
were described (as well as associative ones).
The main result  of~\cite{12} was generalized in~\cite{1}
for finite pseudoalgebras.

Some features of the structure theory
of conformal algebras (and their representations)
of infinite type were also considered
in a series of works
(see
\cite{8,9,13,21,23,24,28,29}).
One of the most urgent problems in this field
is to describe the structure of conformal algebras with faithful
irreducible representation of finite type (these algebras could
be of infinite type themselves). In~\cite{21}
and~\cite{8},
the conjectures on the structure of such algebras (associative
and Lie) were stated. The papers
\cite{8,13,29}
contain confirmations of these conjectures under some additional
conditions.

Another problem is to classify simple and semisimple conformal
algebras
of linear growth (i.e., of Gel'fand--Kirillov
dimension one)~\cite{23}.
In the papers
\cite{23,24,28,29},
this problem was solved for
finitely generated associative conformal algebras
which contain a unit
(\cite{23,24}),
or at least an idempotent
(\cite{28,29}).
The objects appearing from the consideration of
conformal algebras with faithful irreducible
representation of finite type are similar to
those examples of conformal algebras stated in these papers.

A combinatorial aspect of the theory of conformal algebras
requires the notion of a free conformal (and vertex) algebra.
Actually, it was stated in~\cite{7} what is
to be the free vertex algebra.
Studying of free conformal algebras was initiated
in~\cite{25}.
The basic notions of the combinatorial theory of associative
conformal algebras
(Gr\"obner--Shirshov bases, Composition-Diamond lemma)
were considered in
\cite{4,5,6}.

The class of conformal algebras coincides with the class of
pseudoalgebras \cite{1,2} over the polynomial Hopf algebra $\Bbbk[D]$,
where $D$ is just a formal variable,
$\Bbbk$ is a field of zero characteristic.
In general, a pseudoalgebra over a cocommutative Hopf algebra $H$
is just an algebra in the pseudotensor category associated
with~$H$ \cite{2}.
Therefore, the notion of a pseudoalgebra provides
the common point of
view to the theory of ordinary algebras ($H=\Bbbk$) and to
the theory
of conformal algebras ($H=\Bbbk[D]$).
In particular, for every finitely generated $H$-module
$V$ one can define
a pseudoalgebra structure on the set
of all $H$-linear
(with respect to the second tensor factor)
 maps $V\to H\otimes H\otimes_H V$~\cite[Proposition~10.1]{1}.
The pseudoalgebra obtained is denoted by $\Cend V$
(conformal endomorphisms).
This is a direct
analogue of the ordinary algebra $\End U$
of all linear transformations
of a vector space~$U$.

From now on, let
$\Bbbk$ be an algebraically closed field
of zero characteristic, and let $H=\Bbbk[D]$ be the polynomial algebra
in one variable~$D$. By  $\partial_x$ we denote  the operator
of formal  derivation
with respect to a variable~$x$.
We will also use the notation $x^{(n)} = x^n/n!$, $n\in \Zset_+$,
where $\Zset_{+}$ means  the set of non-negative integers.

Our purpose is to develop the structure theory of
associative conformal algebras with faithful representation
of finite type
(as a corollary, we obtain the classification
of simple and semisimple
associative conformal algebras of finite type~\cite{21}).
The main point of this theory is an analogue
of the classical Burnside theorem
on irreducible subalgebras of the matrix algebra
$M_N(\Bbbk)$.
The last statement says that, if a subalgebra
$S\subseteq M_N(\Bbbk)$ acts irreducibly on $\Bbbk^N$
(i.e., there are no non-trivial submodules),
 then $S=M_N(\Bbbk)$.
A similar
  question for conformal algebras was raised by V.~Kac~\cite{21}.
In this context, the ``classical" objects
$\Bbbk$,
$\Bbbk^N$,
$M_N(\Bbbk)$
should be replaced with
$H$,
$V_N$,
$\Cend_N$,
respectively.
Here
$V_N$ is the free $N$-generated $H$-module,
$\Cend_N=\Cend V_N$.
The associative conformal
algebra $\Cend_N$ (and its adjoint Lie conformal algebra $\gc_N$)
 plays the same role in the theory of conformal algebras as
$M_N(\Bbbk)$ (respectively, $\Gll_N(\Bbbk)$) does in the theory
of ordinary algebras.

The paper~\cite{8}
contains a systematical investigation of the algebra
$\Cend_N$: the structure of its left/right ideals and their
modules of finite type, classification
of (anti-) automorphisms and involutions, etc.
Since the algebra $\Cend_N$ is the main object of our study,
we will also state some of these results.

In the paper \cite{21} the following statement
was  conjectured:
{\sl
every irreducible
conformal subalgebra $C\subseteq \Cend_N$ is either
a left ideal of $\Cend_N$
or a conjugate of the current subalgebra}
(in the last case, $C$ is of
finite type).
In \cite{8},
the second part of this conjecture (i.e., finite type case)
was derived
  from
the well-developed representation theory of Lie
conformal algebras of finite type~\cite{12}.
Also, the conjecture was proved in~\cite{8}
 for $N=1$ as well as for unital conformal subalgebras.

We prove this conjecture in general (Theorem \ref{thm 5.5.1}).
Moreover,
although it follows from
\cite{12,21} that
any simple associative conformal algebra  of finite type
is isomorphic to the current conformal algebra over $M_n(\Bbbk)$,
we focus on the independent proof of this fact.

The paper is organized as follows.
In Section~\ref{sec2} we state the definitions of conformal algebras
and their representations.
Also, we write down the general construction of an associative
conformal algebra:
given an (ordinary) associative conformal algebra $A$ with
a locally nilpotent
derivation~$\partial$,
the free $H$-module $H\otimes A$ could be endowed with the
structure of an associative conformal algebra \cite{1,21,23,24}
denoted by $\Diff A$.

In Section~\ref{sec3}, we introduce the conformal algebra
$\Cend_N=\Cend V_N$
($N$ is a positive integer) which is the main object of our study.
This conformal algebra is isomorphic to
$\Diff M_N(\Bbbk[v])$, where $v$ is a formal variable,
and $\partial = \partial_v$ is the usual derivation
with respect to~$v$.
Therefore, $\Cend_N$ could be identified with $M_N(\Bbbk[D,v])$.

As a main tool of our investigation
we use a correspondence between conformal subalgebras
of $\Cend_N$
and subalgebras of $M_N(W)$, where
$W=\Bbbk \langle p,q \mid qp-pq=1\rangle $
is the 1st Weyl algebra.
This correspondence is provided by a construction
called operator algebra,
which is very close to  the  annihilation algebra
\cite{1,21}. Namely, for any conformal subalgebra
$C\subseteq \Cend_N $ we consider the ordinary algebra
 $S(C)\subseteq M_N(W)$ which consists of linear operators on~$V_N$
corresponding to all elements of~$C$.  If $C$ is a left (right)
ideal then so is~$S(C)$. Moreover, $C$ has a canonical structure
of a (left) $S(C)$-module.

In Section~\ref{sec4}, we prove some properties of
operator algebras. First,
we show that for any irreducible conformal subalgebra
$C\subseteq \Cend_N$ its ``expanded" operator algebra $\Bbbk[p]S(C)$
is a dense subalgebra of $M_N(W)$ with respect to the finite
topology~\cite{19}.
The rest of the section is devoted to the properties of ``small"
operator algebras
(i.e., those of linear growth)
corresponding to conformal subalgebras of finite type.
Roughly speaking, we show that if $S=S(C)\subset M_N(W)$
is an operator algebra
of linear growth corresponding to an irreducible
conformal subalgebra $C\subseteq \Cend_N$,
then $C$ is a conjugate of the current
subalgebra of $\Cend_N$.
In this way, we obtain an independent description
of irreducible conformal subalgebras of finite type
\cite[Theorem~5.2]{8}.

Section \ref{sec5} contains the main results of the paper.
The basic point is to describe all irreducible conformal
subalgebras of $\Cend_N$.
First, we show that if $C\subseteq \Cend_N$
is irreducible, then
\begin{equation}
\sum_{n\ge 0}v^n C = \Cend_{N,Q},
                               \label{R1}
\end{equation}
where $\Cend_{N,Q}=M_N(\Bbbk[D,v])Q(-D+v)$,
$Q=Q(v)$ is a nondegenerate
matrix with polynomial entries.

It is left to consider three cases:

1) the sum (\ref{R1}) is direct;

2) $C\cap vC \ne 0$;

3) the sum (\ref{R1}) is non-direct, but $C\cap vC =0$.

\noindent
The first case is considered in Section~\ref{subsec5.2}.
Using the results
of Sections~\ref{subsec4.2} and~\ref{subsec4.3},
we show that $C$ is a conjugate of
the current subalgebra of $\Cend_N$.
In the second case, we use the method proposed in \cite[Section~2]{8}
for $N=1$: in Section~\ref{subsec5.3}
we reduce the computation to the one-dimensional
case and obtain $C=\Cend_{N,Q}$.
The third case is impossible, as we show in Section \ref{subsec5.4}.

In Section~\ref{subsec5.5} we complete the classification of
irreducible
subalgebras (Theorem \ref{thm 5.5.1}) and proceed with the
structure theory of associative conformal algebras with faithful
representation of finite type: we describe simple and semisimple
ones. Moreover, we show that any conformal algebra of this kind
has a maximal nilpotent ideal. As a corollary, we obtain the
classification theorem for associative conformal algebras of
finite type~\cite{21}.

\section{Conformal algebras and their representations}\label{sec2}

In this section, we state necessary definitions and constructions
from the theory of conformal algebras.

\subsection{Definition and examples of conformal
algebras}\label{subsec2.1}

\begin{defn}[\cite{20}]\label{defn 2.1.1}
A vector space $C$ endowed with a linear map
$D$ and with a family of bilinear operations
$\oo{n}$, $n\in \Zset_+$, is said to be a {\em conformal algebra},
 if
\begin{eqnarray*}
\hbox{\rm (C1)} && \quad a\oo{n}b=0,\quad
\hbox{for $n$ sufficiently large},\quad
   a,b\in C; \\
\hbox{\rm (C2)} && \quad Da\oo{n}b=-n a\oo{n-1} b; \\
\hbox{\rm (C3)} && \quad a\oo{n} Db = D(a\oo{n}b)+na\oo{n-1}b.
\end{eqnarray*}
\end{defn}
The axiom (C1) allows to define so called {\em locality function}
${\mathcal N}:C\times C \to \Zset_+$. Namely,
\begin{equation}\label{defn LocFunct}
{\mathcal N}(a,b) = \min\{N\in \Zset_+ \mid a\oo{n}b = 0\
  \mbox{for all}\ n\ge N\},
\quad a,b \in C.
\end{equation}

Every conformal algebra could be considered as a unital left
module over $H=\Bbbk[D]$ endowed with a family of
sesquilinear (i.e., satisfying (C2), (C3)) products
$\oo{n}$, $n\in \Zset_+$,
such that the locality axiom (C1) holds.

For any pair of $D$-invariant subspaces
($H$-submodules) $X,Y\subseteq C$
the space $X\oo{\omega} Y =\sum\limits_{n\in \Zset_+}X\oo{n}Y$
is also an $H$-submodule.

An ideal of a conformal algebra $C$ is
an $H$-submodule $I\subseteq C$
which is
closed under all multiplications by $C$, i.e.,
$C\oo{\omega} I, I\oo{\omega} C \subseteq I$.
Left and right ideals are defined analogously.

A conformal algebra $C$
is simple, if
$C\oo{\omega} C\ne 0$
and there are no non-trivial ideals.

For any ideal $I$ of a conformal algebra $C$ one can define
descending sequences of ideals
$\{I^{(n)}\}_{n=1}^\infty $
and
$\{I^{n}\}_{n=1}^\infty $
by the usual rule:
$I^{(1)}=I^1=I$,
$I^{(n)}=I^{(n-1)}\oo{\omega} I^{(n-1)}$,
$I^n = \sum_{s=1}^{n-1}I^s\oo{\omega } I^{n-s}$, $n>1$.
An ideal $I$ is said to be
solvable (resp., nilpotent),
if $I^{(n)}=0$ (resp., $I^n=0$) for sufficiently large~$n$.
If a conformal algebra $C$ has no non-zero solvable ideals
then $C$ is called semisimple.

If a conformal algebra $C$ is finitely generated
as an $H$-module, then it is said to be
 {\em of finite type\/}
(or {\em finite\/} conformal algebra).

A homomorphism $\phi : C_1\to C_2$ of conformal algebras  is
a $D$-invariant linear map
such that
$\phi(a\oo{n} b) = \phi(a)\oo{n}\phi (b)$,
$a,b\in C_1$, $n\in \Zset_+$.

For any conformal algebra
$C$ in the sense of
Definition~\ref{defn 2.1.1} there exists an ordinary algebra
$A$ such that $C$ lies in
the space of formal distributions
$A [[z,z^{-1}]]$, where $D=\partial_z$ and
the $\oo{n}$-products are given by
\begin{equation}\label{2.1.1}
a(z)\oo{n} b(z) = \Ress_{w=0} a(w)b(z)(w-z)^n,
\quad n\in \Zset_+.
\end{equation}
Such an algebra $A$ is not unique, but there exists a universal one
denoted by $\Coeff C$.
Namely, for any algebra $A$
such that $A[[z,z^{-1}]]$ contains $C$ as above, there exists a
homomorphism $\Coeff C \to A$
such that the natural expansion
$\Coeff C[[z,z^{-1}]]\to A[[z,z^{-1}]]$
acts on $C$ as an identity \cite{21,25}.
The algebra $\Coeff C$ is called the
{\em coefficient algebra\/} of~$C$.

There is a correspondence between identities on $\Coeff C$
and systems of conformal identities on~$C$.
In particular, $\Coeff C$
is associative if and only if $C$ satisfies
\begin{equation}\label{2.1.2}
(a\oo{n} b)\oo{m} c = \sum\limits_{s\ge 0} (-1)^s
 \binom{n}{s}
a\oo{n-s} (b\oo{m+s} c)
\end{equation}
or
\begin{equation}\label{2.1.3}
a\oo{n} (b\oo{m} c) = \sum\limits_{s\ge 0}
 \binom {n}{s} (a\oo{n-s} b)\oo{m+s} c,
\end{equation}
for all  $n,m\in \Zset_+$.
The systems of relations (\ref{2.1.2}) and (\ref{2.1.3})
are equivalent.
Also, $\Coeff C$ is a Lie algebra if and only if
$C$ satisfies
\begin{eqnarray}
& a\oo{n} b = -\sum\limits_{s\ge 0} (-1)^{n+s}D^{(s)}(b\oo{n+s} a),
   \label{2.1.4}  \\
& a\oo{n}(b\oo{m} c) - b\oo{m}(a\oo{n} c) =
 \sum\limits_{s\ge 0}
  \binom{n}{s}
   (a\oo{n-s} b)\oo{m+s} c
      \label{2.1.5}
\end{eqnarray}
for all $n,m\in \Zset_+$.

A conformal algebra $C$ is called associative, if   $\Coeff C$
is associative, i.e., if $C$ satisfies (\ref{2.1.2})
or (\ref{2.1.3}).
Analogously, $C$ is called Lie conformal algebra,
if $\Coeff C$ is a Lie algebra, i.e., if $C$ satisfies
(\ref{2.1.4}) and (\ref{2.1.5}).
It is easy to note that the notion of solvability coincides
with nilpotency in the case of associative conformal algebras.

\begin{prop}[see, e.g., \cite{21}]\label{prop 2.1.2}
Let $C$ be an associative conformal algebra. Then the
same $H$-module $C$ endowed with new operations
\begin{equation}\label{2.1.6}
[a\oo{n}b]= a\oo{n} b-\sum\limits_{s\ge 0}
   (-1)^{n+s}D^{(s)}(b\oo{n+s} a),
\quad a,b\in C,\ n\in \Zset_+,
\end{equation}
is a Lie conformal algebra denoted by $C^{(-)}$.
\quad \qed
\end{prop}

\begin{exmp}\label{exmp 2.1.3}
Let $\mathcal A$ be an (associative or Lie) algebra. Then
 formal distributions of the form
\[
a(z)=\sum\limits_{n\in \Zset} at^n z^{-n-1} \in
   {\mathcal A}[t,t^{-1}][[z,z^{-1}]],
   \quad a\in {\mathcal A},
\]
together with all their derivatives
span a conformal algebra in ${\mathcal A}[t,t^{-1}][[z,z^{-1}]]$
called the {\em current conformal algebra\/}
$\Curr {\mathcal A}$.
The operations (\ref{2.1.1}) are given by
\[
a(z)\oo{n}b(z) = \begin{cases}
(ab)(z), & n=0; \\
                       0, & n>0.
                       \end{cases}
\]
The conformal algebra $\Curr {\mathcal A}$ is associative or Lie
if and only if
so is~${\mathcal A}$.
\end{exmp}

\begin{exmp}\label{exmp 2.1.4}
Consider ${\mathcal A}=\Bbbk\langle t,t^{-1},\partial \mid
\partial t-t\partial =1 \rangle$.
Then formal distributions of the form
\[
v_m (z) = \sum\limits_{n\in \Zset} t^n\partial^m z^{-n-1}
\in {\mathcal A}[[z,z^{-1}]],
\quad m\in \Zset_+,
\]
together with all their derivatives span a conformal algebra in
${\mathcal A}[[z,z^{-1}]]$ called  the
{\em Weyl conformal algebra}.
\end{exmp}

\begin{exmp}\label{exmp 2.1.5}
Consider the element $L=v_1$ from the previous example. It is easy
to compute that the one-generated $H$-submodule $H\otimes \Bbbk L$ of
the Weyl conformal algebra is closed under
  operations (\ref{2.1.6}):
\[
[L\oo{0}L] = -DL, \quad [L\oo{1} L ] = -2L,
\quad
[L\oo{n} L] = 0,\ n\ge 2.
\]
The Lie conformal algebra obtained is called the {\em Virasoro
conformal algebra}.
\end{exmp}


\subsection{Representations of associative conformal
algebras}\label{subsec2.2}

Remind that $H = \Bbbk [D]$ is a Hopf algebra,
so  $H\otimes H$  could be considered as the
outer product  of regular right $H$-modules
(the action is given by
$(f\otimes g)\cdot D = fD\otimes g + f\otimes gD$).
If $\{h_i\mid i\in I\}$ is a linear basis of $H$, then
$\{h_i\otimes 1 \mid i\in I\}$ is an $H$-basis of the right $H$-module
$H\otimes H$ (see, e.g., \cite[Lemma~2.3]{1}).

\begin{defn}[\cite{1}]\label{defn ConfEnd}
Let $V$ be a unital left $H$-module.
A linear map
\[
a : V \to (H\otimes H)\otimes_H V
\]
is said to be a {\em conformal endomorphism\/} of $V$,
if
$a(fv)= ((1\otimes f)\otimes_H 1)a(v)$
for any $f\in H$, $v\in V$.
\end{defn}

Denote the set of all conformal endomorphisms of $V$ by $\Cend V$.
For any $a\in \Cend V$ there exists a uniquely
defined sequence $\{a_n\}_{n= 0}^\infty$
of (ordinary) $\Bbbk $-linear endomorphisms of $V$ such that:
\begin{eqnarray}
& a_n(v)=0 \quad  \mbox{for $n$ sufficiently large},\quad v\in V;
    \label{CV1}  \\
& [a_n, D] = na_{n-1},  \quad n\in \Zset_+; \label{CV2-3}\\
& a(v) = \sum\limits_{n\ge 0} ((-D)^{(n)}\otimes 1) \otimes _H a_n(v),
 \quad v\in V.
                                          \label{CV}
\end{eqnarray}
Conversely, any sequence $\{a_n\}_{n= 0}^\infty \subset \End_{\Bbbk} V$
satisfying (\ref{CV1}) and (\ref{CV2-3})
defines a conformal endomorphism $a\in \Cend V$ via (\ref{CV}).

For any $a,b\in \Cend V$ define the sequences
$\{(Da)_n \}_{n= 0}^\infty$
and $\{(a\oo{m} b)_n\}_{n= 0}^\infty$, $m\in \Zset_+$,
by the following rules (see, e.g., \cite[Section~3]{12}):
\begin{eqnarray}
&   (Da)_n = [D, a_n] = -n a_{n-1}, \label{DCend} \\
&   (a\oo{m}b)_n = \sum\limits_{s\ge 0} (-1)^s \binom{m}{s}
   a_{m-s} b_{n+s}.
 \label{nCend}
\end{eqnarray}
It is clear that these sequences satisfy
(\ref{CV1}) and (\ref{CV2-3}),
so the conformal endomorphisms $Da$, $a\oo{m} b$
($m\in \Zset_+$) are defined.

For an arbitrary $H$-module $V$, the vector space $\Cend V$
endowed with the operations $D$ and $\oo{m}$ ($m\in \Zset_+$)
given by (\ref{DCend}) and (\ref{nCend})
satisfies (C2) and (C3). The locality property (C1)
does not hold, in general.
However,  if $V$ is a finitely generated $H$-module, then (C1) holds,
so $\Cend V$ is an associative conformal algebra.

Let $C$ be a conformal algebra, and let $V$ be a left unital
$H$-module.
If $\varrho : C \to \Cend V$ is a linear map
such that
\[
\varrho(Da)=D\varrho(a), \quad
\varrho(a\oo{n} b) = \varrho(a)\oo{n} \varrho(b),
\quad a,b\in C, \ n\in \Zset_+,
\]
then $\overline{C}=\varrho(C)\subseteq \Cend V$ is a conformal algebra.
By abuse of terminology, we will say that
$\varrho $ is a homomorphism of conformal algebras
(although $\Cend V$ is not necessarily a conformal algebra).

\begin{defn}[\cite{1,10,20,21}]\label{defn 2.4.1}
A {\em representation\/} of an associative conformal algebra $C$
on a left unital $H$-module $V$
is a homomorphism $\varrho : C\to \Cend V$ of
conformal algebras.
If $C$ has a representation on $V$, then  $V$ is said
to be a module over conformal algebra~$C$
(or $C$-module, if the context is clear).
\end{defn}

The last definition is equivalent to the following one \cite{10}:
an $H$-module $V$ is a module over an associative conformal algebra
$C$, if there is a family of $\oo{n}$-products
\begin{equation}\label{n-action}
\oo{n} : C\otimes V \to V
\end{equation}
such that
\begin{eqnarray*}
\hbox{\rm (M1)} && \quad a\oo{n}v=0,\quad
 \hbox{for $n$ sufficiently large},\quad
   a\in C,\ v\in V; \\
\hbox{\rm (M2)} && \quad Da\oo{n}v=-n a\oo{n-1} v; \\
\hbox{\rm (M3)} && \quad a\oo{n} Dv = D(a\oo{n}v)+na\oo{n-1}v; \\
\hbox{\rm (M4)} && \quad a\oo{n}(b\oo{m} v) =
  \sum\limits_{s\ge 0} \binom{n}{s} (a\oo{n-s} b) \oo{m+s} v.
\end{eqnarray*}

A representation $\varrho: C\to \Cend V$ is called faithful,
if $\varrho $ is injective.
If  $V$  is  a  finitely generated  $H$-module,
then $\varrho $ is said to be a representation of finite type
(or finite representation).

In the case of ordinary algebras, every algebra which has
a faithful finite-dimensi\-onal representation is
finite-dimensional itself. It is not the case for
conformal algebras: there exist infinite conformal algebras
with faithful finite representation.


\subsection{Differential conformal algebras}\label{subsec2.3}

Let $\mathcal A$ be an associative algebra with a locally
nilpotent derivation $\partial $. We consider the free $H$-module
$H\otimes {\mathcal A}$ as an associative  conformal algebra.
Notice that by (C2), (C3) it is sufficient to define the family of
$\oo{n}$-products between elements of the form $1\otimes a$, $a\in
{\mathcal A}$. Hereinafter we identify $1\otimes a$ with $a\in
{\mathcal A}$.

\begin{prop}[see, e.g., \cite{1,22,23}]\label{prop 2.3.1}
{\em (i)}
The family of $\oo{n}$-products
\begin{equation}\label{2.3.1}
a\oo{n}b = a \partial^n(b) ,
\quad
a,b\in {\mathcal A},
\end{equation}
defines an associative conformal algebra structure on
$H\otimes {\mathcal A}$.
The conformal algebra obtained is denoted by $\Diff {\mathcal A}$.

{\em (ii)}
The same is true for another family of operations:
\begin{equation}\label{2.3.2}
(a\oo{(n)} b)= \sum\limits_{s\ge 0} D^{(s)}\otimes \partial^{n+s}(a)b,
\quad
a,b\in {\mathcal A}.
\end{equation}
The conformal algebra obtained is denoted by
$\Diff^{\circ} {\mathcal A}$.
 \quad  \qed
\end{prop}

\begin{prop}\label{prop 2.3.2}
$\Diff {\mathcal A}\simeq \Diff^\circ {\mathcal A}$.
\end{prop}

\begin{proof}
The isomorphism
$\varphi : \Diff {\mathcal A} \to  \Diff^\circ {\mathcal A}$
is given by
\[
h\otimes a \mapsto \sum\limits_{s\ge 0}D^{(s)}h \otimes \partial^s(a);
\]
the inverse
$\varphi^{-1}: \Diff^\circ {\mathcal A}\to \Diff {\mathcal A}$
could be defined via
\[
h\otimes a    \mapsto \sum\limits_{s\ge 0}
 (-D)^{(s)}h \otimes \partial^s(a).
\]
Direct computation shows that
$\varphi (a\oo{n} b)=\varphi(a) \oo{(n)} \varphi(b)$.
\end{proof}

In \cite{23,24}, an associative conformal algebra obtained by
either of the constructions (\ref{2.3.1}), (\ref{2.3.2})
is called {\em differential}.

\begin{exmp}\label{exmp 2.3.3}
Let ${\mathcal A}$ be an (associative) algebra,
and let $\partial$ be the zero derivation:
$\partial({\mathcal A})= 0$. Then
$\Diff {\mathcal A}=\Diff^\circ {\mathcal A}$
is exactly the current conformal algebra
$\Curr {\mathcal A}$ (Example~\ref{exmp 2.1.3}).
\end{exmp}

\begin{exmp}\label{exmp 2.3.4}
Consider
${\mathcal A}=\Bbbk[v]$,
 and let $\partial =\partial_v$ be the usual
derivation. Then
the conformal algebras
$\Diff {\mathcal A}$ and $\Diff^\circ {\mathcal A}$
are isomorphic to the Weyl conformal algebra
(Example~\ref{exmp 2.1.4}).
\end{exmp}

\begin{exmp}\label{exmp 2.3.5}
Let $X$ be a set of symbols (generators), $v\notin X$, and let
${\mathcal A}=\Bbbk\langle X, v \rangle$
be the free associative algebra generated by $X\cup\{v\}$
endowed with the derivation $\partial =\partial_v$.
Consider a conformal subalgebra of $\Diff {\mathcal A}$
 generated by the set
$\{v^{(N-1)}x \mid x\in X\}$, $N\ge 1$.
This subalgebra is isomorphic to
the {\em free associative conformal algebra\/}
generated by $X$ with respect to the constant
locality function ${\mathcal N}(x,y)\equiv N$,
$x,y\in X$
(see \cite{4,25}).
\end{exmp}

The last example shows that every finitely generated
associative conformal
algebra is a homomorphic image of a subalgebra of a
differential conformal algebra.

It was conjectured in \cite{24}
that every annihilation-free (see \cite{23} for rigorous explanation)
associative conformal
algebra is just a subalgebra of some differential conformal
algebra. Example \ref{exmp 2.3.5}
confirms this conjecture for free associative
conformal algebras.


\section{Conformal algebra  $\Cend_N$}\label{sec3}

From now on, we will use the term
``conformal algebra" for associative conformal algebra,
unless stated otherwise.

\subsection{Construction}\label{subsec3.1}

Let us introduce the main object of our study: the conformal
algebra $\Cend_N =\Cend V_N$, where $V_N$ is the free
$N$-generated $H$-module.
We will use the following presentation of $\Cend_N$.

Consider the algebra $M_N(\Bbbk[v])$ of all
$N\times N$-matrices with polynomial entries
 endowed  with the derivation $\partial _v$
with respect to the variable~$v$.
Denote the conformal algebras $\Diff M_N(\Bbbk[v])$
and $\Diff ^\circ M_N(\Bbbk[v])$ by
${\mathfrak A}_N$ and ${\mathfrak A}_N^{\circ}$, respectively.
Also, there is the usual multiplication on $H\otimes M_N(\Bbbk[v])$
as on $M_N(\Bbbk[D,v])$, and
the following property of ${\mathfrak A}_N$ is clear:
\begin{equation}\label{3.1.1}
v(a\oo{n} b) = va\oo{n} b,
\quad
a\oo{n}vb = va\oo{n} b + n a\oo{n-1}b.
\end{equation}

Let us fix an $H$-basis
$\{ e_1,\dots,e_N \}$
of $V_N$, and
let $U\subset V_N$ be the finite-dimensional subspace spanned
by this basis.
Since $U\simeq \Bbbk^N $, the free $H$-module
$V_N$ could be identified with $H\otimes \Bbbk^N = H^N$.
For any element
$a=\sum_{s\ge 0}(-D)^{(s)}\otimes A_s(v)\in {\mathfrak A}_N$,
$A_s(v)\in M_N(\Bbbk[v])$,
consider the following sequence of
$\Bbbk$-linear endomorphisms of $V_N$:
\begin{equation}\label{3.1.2}
a_n : u \mapsto \sum\limits_{s\ge 0}\binom{n}{s}
A_s(D) \partial_D^{n-s}(u),
\quad
 u\in V_N, \ n\in \Zset_+.
\end{equation}
It is easy to see that this sequence satisfies
(\ref{CV1}) and (\ref{CV2-3}),
so it defines a conformal endomorphism $\tilde a\in \Cend_N$.

\begin{prop}[\cite{1,12,20,23}]\label{prop 3.1.1}
The map $a\mapsto \tilde a$
is an isomorphism of conformal
algebras ${\mathfrak A}_N$ and $\Cend_N$.
\quad \qed
\end{prop}

\begin{rem}\label{rem 3.1.2}
It follows from Proposition \ref{prop 2.3.2} that
\[
{\mathfrak A}_N \simeq {\mathfrak A}_N^{\circ }\simeq \Cend_N.
\]
The isomorphism $\varphi $ from Proposition \ref{prop 2.3.2}
corresponds to the map ${\mathfrak A}_N \to {\mathfrak A}_N^\circ$
given by $D\mapsto D$, $v\mapsto D+v$.
In \cite{8,21},
the second construction (\ref{2.3.2}) was used for $\Cend_N$.
We will preferably follow \cite{23}
and use the $\oo{n}$-products (\ref{2.3.1}) for $\Cend_N$.
\end{rem}

From now on, identify
$\Cend_N$ with~${\mathfrak A}_N = H\otimes M_N(\Bbbk[v])$,
where $A(v)\oo{n} B(v) = A(v)\partial_v^n(B(v))$.

\begin{defn}\label{defn 3.1.3}
A conformal subalgebra $C\subseteq \Cend_N$ is called
{\em irreducible},
if there are no non-trivial $C$-submodules of~$V_N$.
\end{defn}


\subsection{Operator Algebras}\label{subsec3.2}

Let us fix an $H$-linear basis of $V_N$, and identify
$V_N$ with $H\otimes \Bbbk^N$ as above.
For every $a\in \Cend_N$ and for every $n\in \Zset_+$
consider the  $\Bbbk$-linear map
\[
a(n)=a_n : V_N \to V_N, \quad a(n): u\mapsto a\oo{n}u,\ u\in V_N,
\]
given by (\ref{3.1.2}).
Introduce a new variable $p$ instead of $D$ in this
context, and  set  $q = \partial_p$.
Thus, we get
$a(n)\in M_N(W)$, where $W=\Bbbk\langle p,q\mid qp-pq=1  \rangle$
is the (1st) Weyl algebra.
 The axiom (C2) implies
 $(Da)(n)=-na(n-1)=-\partial_q a(n)$,
 where $\partial_q$ is the partial derivation with respect to~$q$.

In fact, the set
$S(\Cend_N)=\{a(n)\mid a\in \Cend_N, \, n\in \Zset_+\}$
is a linear space with
multiplication (composition) given by
\begin{equation}\label{3.2.1}
a(n)b(m)=\sum\limits_{s\ge 0}\binom{n}{s}
   (a\oo{n-s} b)(m+s).
\end{equation}
It is easy to note that
$S(\Cend_N) = M_N(W)$.

For every $H$-submodule $X$ of $\Cend_N$ one may consider
the subspace
\[
S(X)=\{a(n) \mid a\in X, \, n\in \Zset_+\}\subseteq S(\Cend_N).
\]
If $C$ is a conformal subalgebra of $\Cend_N$,
then $S(C)$
is a $\partial_q$-invariant
 subalgebra of $M_N(W)$. The same is true for left and right ideals.

\begin{defn}\label{defn 3.2.1}
Let $C \subseteq \Cend_N$ be a conformal subalgebra.
The subalgebra $S(C)$ of $M_N(W)$
is called the {\em operator algebra\/} of~$C$.
\end{defn}

Let us consider the topology
(called {\em $q$-topology}, for short) on $M_N(W)$
defined by the sequence of the
left ideals generated by $q^n$, $n\ge 0$:
\[
M_N(W)\supset M_N(W)q \supset \dots \supset
M_N(W)q^n \supset \dots \supset 0.
\]
It is clear that this topology is equivalent to the
finite topology~\cite{19}
defined by the canonical action of
$M_N(W)$ on~$V_N$.

\begin{defn}\label{defn 3.2.2}
A sequence $\{a_n\}_{n=0}^\infty $,
$a_n\in M_N(W)$,
is called {\em differential},
if $\partial_q(a_n)=n a_{n-1}$
(we mean $\partial_q(a_0)=0$) and
$\lim\limits_{n\to \infty} a_n = 0$
(in the sense of $q$-topology).
\end{defn}

\begin{lem}\label{lem 3.2.3}
Let $\{a_n\}_{n=0}^\infty $ be a differential sequence.
Then there exist a finite number of matrices
$A_s\in M_N(\Bbbk[p])$, $s=0,\dots ,m$,
such that
\[
a_n=\sum\limits_{s= 0}^m \binom{n}{s} A_s q^{n-s}
\]
for all $n\in \Zset_+$.
\end{lem}

\begin{proof}
By the definition of a differential sequence, there exists an
integer $m\ge 0$ such that
$a_n\in M_N(W)q$
for all $n> m$.
Then we may represent the element $a_m$ in the following form:
\[
a_m = \sum\limits_{s=0}^m \binom{m}{s} A_sq^{m-s},
\]
where
$A_0,\dots, A_m \in M_N(\Bbbk[p])$.
For every $n\ge 0$ we may write
\[
a_n=\sum\limits_{s=0}^m \binom{n}{s} A_s q^{n-s}
 + \sum\limits_{k\ge m+1} \binom{n}{k} B_k q^{n-k}.
\]
The conditions $a_n\in M_N(W)q$
(for $n>m$)
and $\partial_q a_n = na_{n-1}$
  imply $B_k=0$ for all $k>m$.
\end{proof}

\begin{prop}\label{prop 3.2.4}
A subalgebra $S\subseteq M_N(W)$ is an
operator algebra of some conformal subalgebra
$C\subseteq \Cend_N$
if and only if every element of $S$ lies in
a differential sequence
of elements from~$S$.
\end{prop}

\begin{proof}
If $S=S(C)$ for some $C\subseteq \Cend_N$, then every
$a(m)\in S$  ($a\in C$, $m\ge 0$)
lies in the  differential sequence
$\{a(n)\}_{n=0}^\infty $.

Let us prove the converse.
By Lemma \ref{lem 3.2.3}, every differential
sequence $\{a_n\}_{n=0}^\infty $
corresponds to a finite family of
matrices $A_s\in M_N(\Bbbk [p])$. Consider the elements
\[
a = \sum\limits_{s\ge 0} (-D)^{(s)}\otimes A_s(v) \in \Cend_N
\]
corresponding to all differential sequences in~$S$.
It is clear that these elements form a conformal
subalgebra~$C$
of $\Cend_N$.
Indeed, let  $a,b\in \Cend_N$ correspond to
$\{a_n\}_{n=0}^\infty $ and
$\{b_n\}_{n=0}^\infty $, respectively.
Then $Da$ corresponds to $\{-\partial_q(a_n)\}_{n=0}^\infty $,
and $a\oo{m} b$ corresponds to $\{c_n\}_{n=0}^\infty $,
where
$c_n = \sum\limits_{s\ge 0 } (-1)^s \binom {m}{s} a_{m-s}b_{n+s}$
(it is easy to see that these $c_n$'s form a differential
sequence).
It follows from the construction that $S(C)=S$.
\end{proof}

\begin{defn}\label{defn 3.2.5}
Let $S\subseteq M_N(W)$ be a subalgebra satisfying the conditions
of Proposition \ref{prop 3.2.4}.
Then the conformal algebra $C$ constructed
is called the {\em reconstruction\/}
of $S$ (c.f. \cite[Section~11]{1}).
We will denote it by
${\mathcal R}(S)$.
\end{defn}

In particular, every operator algebra $S=S(C)\subseteq M_N(W)$
is $\partial_q$-invariant. Hence, the subspace
$\Bbbk[p]S = S+pS+p^2S+\dots $
is also a subalgebra of $M_N(W)$.

\begin{prop}\label{prop 3.2.6}
A conformal subalgebra $C\subseteq \Cend_N$
is irreducible if
and only if
the subalgebra $S_1=\Bbbk[p]S(C)$
acts irreducibly on~$V_N$.
\end{prop}

\begin{proof}
Let $C$ be an irreducible conformal subalgebra.
For every $0\ne u \in V_N$
the following $H$-module is also a $C$-module:
\[
V_u = H\{a\oo{n} u \mid a \in C, \, n \in \Zset_+ \}.
\]
If $V_u=0$, then $\{u\in V_N \mid V_u=0\}$
is a non-trivial $C$-submodule.
Hence, $V_u=V_N$.
Then for every $w\in V_N$ there exists
a finite family
$\{a_i, n_i, f_i\}$,
$a_i\in C$, $n_i\in \Zset_+$, $f_i\in H$,
such that
$w=\sum\limits_i f_i (a_i \oo{n_i} u)$.
Then
for the $S$-module $V_N$
we have
$w = \sum\limits_i f_i(p) a_i(n_i) u \in S_1 u$.

The converse statement has a similar proof.
\end{proof}

\begin{prop}\label{prop 3.2.7}
The algebra $M_N(W)$ acts on the vector space $\Cend_N$.
The action is provided by
\begin{equation}\label{3.2.2}
a(n)\cdot b = a\oo{n}b,
\quad
a,b\in \Cend_N,\ n\in \Zset_+.
\end{equation}
\end{prop}

\begin{proof}
The explicit expression for the action is
\[
A(p)q^n \cdot (D^{(s)}\otimes B(v))
=\sum\limits_{t\ge 0} \binom{n}{t}
D^{(s-t)}\otimes A(v)\partial_v^{n-t}B(v).
\]
It is  straightforward to check
the relation (\ref{3.2.2}) and the associativity.
\end{proof}

\begin{cor}\label{cor 3.2.8}
For every conformal subalgebra
$C\subseteq \Cend_{N}$ we have
$S(C)\cdot C \subseteq C$.
\quad \qed
\end{cor}

\begin{cor}\label{cor 3.2.9}
If $C\subseteq \Cend_N$ is a conformal subalgebra,
then $C$ is a left ideal of ${\mathcal R}(S(C))$.
\end{cor}

\begin{proof}
Let $x\in {\mathcal R}(S(C))$ be an element,
corresponding to a differential
sequence $\{b_n\}_{n=0}^\infty $, $b_n\in S(C)$.
Then $x(n)=b_n$
for all $n\ge 0$. For every $a\in C$ we have
$x\oo{n} a = b_n\cdot a \in C$.
\end{proof}


\subsection{Automorphisms of $\Cend_N$ and $M_N(W)$}\label{subsec3.3}

In this section, we consider a relation between automorphisms
of $\Cend_N$ and $M_N(W)$. We will also obtain the full
description of automorphisms of $\Cend_N$, as in \cite{8}.

\begin{prop}\label{prop 3.3.1}
Every automorphism $\Theta $ of the conformal algebra $\Cend_N$
induces an automorphism $\theta $ of $M_N(W)$:
\begin{equation}\label{3.3.1}
\theta : a(n) \mapsto \Theta(a)(n).
\end{equation}
This automorphism is $\partial_q$-invariant and
continuous in the sense of $q$-topology.
\end{prop}

\begin{proof}
First, let us show that the map (\ref{3.3.1}) is well-defined.
Suppose $a(n)=b(m)\in M_N(W)$. Without loss of generality
we may assume $n=m$, so it is sufficient to show
that $a(n)=0$ implies $\Theta(a)(n)=0$.
If $a(n)=0$, then $a = D^{n+1}x$ for some $x\in \Cend_N$,
and $\Theta(a) = D^{n+1}\Theta (x)$.
Hence,
$\Theta(a)(n)=0$ and $\theta $ is well-defined.

It follows from (\ref{3.2.1}) that $\theta $
is an
automorphism of $M_N(W)$. Moreover,
since $\{a(n)\}_{n=0}^\infty$ is a differential sequence,
we may write
$\theta (\partial_q(a(n))) = n \theta (a(n-1))
 =n\Theta(a)(n-1) = \partial_q (\theta (a(n))) $.
Hence,  $\theta $ is $\partial_q$-invariant.

Let us consider the image of $e = 1\otimes \Id_N$
under $\Theta$.
If $\Theta (e) = a\in \Cend_N$, then
$\theta(q^n)=a(n)$. The sequence $\{a(n)\}_{n=0}^\infty $
is differential, and $\theta (M_N(\Bbbk[p]))=M_N(\Bbbk[p])$
since $\theta $ is $\partial_q$-invariant. It is easy
to note that, under these conditions,
$\theta $ has to be continuous in the sense of $q$-topology.
\end{proof}

\begin{exmp}\label{exmp 3.3.2}
Let us consider the map
$\Theta_{\alpha,Q} : \Cend_N \to \Cend_N$
such that
\begin{equation}\label{3.3.2}
\Theta _{\alpha,Q} (a(D,v)) =
 (1\otimes Q^{-1}(v)) a(D,v+\alpha)  Q(-D\otimes 1 + 1\otimes v),
\end{equation}
$\alpha \in \Bbbk$,
$Q,Q^{-1}\in  M_N(\Bbbk[v])$.
It is straightforward to check that
$\Theta_{\alpha,Q}$ is
an automorphism of $\Cend_N$ (see \cite{8}).
The corresponding automorphism $\theta_{\alpha, Q} $ of
$M_N(W)$ acts as follows:
\[
\theta_{\alpha, Q} : a(p,q) \mapsto Q^{-1}(p) a(p+\alpha , q) Q(p).
\]
\end{exmp}

It is shown in \cite{8} that $\Theta_{\alpha, Q}$ exhaust
all automorphisms of $\Cend_N$. Later we will also obtain this
result.

\begin{lem}\label{lem 3.3.3}
Let $\theta $ be a
$\partial_q$-invariant automorphism of $M_N(W)$.
Then
$\theta (p) =p+\alpha$, $\theta(q) = q-B(p)$, where
$\alpha \in \Bbbk$, $B(p)\in M_N(\Bbbk[p])$.
\end{lem}

\begin{proof}
Since $\theta (M_N(\Bbbk[p])) = M_N(\Bbbk[p])$,
we have $\theta(p)\in M_N(\Bbbk[p])$.
The relation
$[p,M_N(\Bbbk[p])]=0$ implies $\theta(p)$ to
be a scalar matrix: $\theta(p)=f(p)\Id_N \equiv f(p)$.
Analogously,
$\partial_q(\theta(q)) = 1$,
so $\theta(q) = q - B(p)$
for some $B(p) \in M_N(\Bbbk[p])$.
Then $\theta([q,p])=[q-B(p), f(p)] = f'(p)=1$,
so, $\deg f=1$.
\end{proof}

\begin{lem}\label{lem 3.3.4}
{\rm (i)}
Let $\theta $ be an automorphism
of $M_N(\Bbbk[p])$ such that $\theta (p)= p$. Then there exists
a nondegenerate matrix $Q\in M_N(\Bbbk[p])$
such that $Q^{-1}\in M_N(\Bbbk[p])$
and
$\theta(X)=Q^{-1}XQ$ for all $X\in M_N(\Bbbk[p])$.

{\rm (ii)}
Let $\theta $ be an automorphism of $M_N(W)$ such that
$\theta (p)= p$.
Then $\theta(q)=q - B(p)$,
where
$B=h(p) - Q^{-1}Q' $
for some polynomial $h(p)\in \Bbbk[p]$ and
invertible matrix
$Q\in M_N(\Bbbk[p])$, $Q' = \partial_p (Q)$.
\end{lem}

\begin{proof}
(i)
Denote $S_0=\theta ^{-1}(M_N(\Bbbk))$. It is easy to
see that
$\Bbbk[p]S_0=M_N(\Bbbk[p])$
and the sum
\begin{equation}\label{3.3.3}
\sum\limits_{n\ge 0}p^nS_0=M_N(\Bbbk[p])
\end{equation}
is direct.

Let us consider the canonical action of $S_0$ on $V_N$.
Relation (\ref{3.3.3}) implies that $S_0u \ne 0$
for every
$0\ne u\in V_N$. So there exist non-zero finite dimensional
$S_0$-submodules of $V_N$,
and we may choose
a non-zero $S_0$-submodule $U_0$ of minimal dimension over~$\Bbbk$.
Since $S_0$ is simple, the $S_0$-module $U_0$
is faithful and $\dim _{\Bbbk}U_0=N$.

Consider $u_1,\dots, u_N\in V_N$ such that
$U_0=\Bbbk u_1 \oplus \dots \oplus \Bbbk u_N$.
Then these vectors are linearly independent over
the field of rational functions
$\Bbbk(p)$. Indeed, let there exist
$f_1(p),\dots, f_N(p)\in \Bbbk(p)$ such that
\begin{equation}\label{3.3.4}
f_1(p)u_1 + \dots + f_N(p)u_N = 0
\end{equation}
(it is sufficient to consider only $f_i(p)\in \Bbbk[p]$).
Then we may choose $a_i\in S_0=\End_{\Bbbk} U_0$
such that
$a_i u_j = \delta_{ij}u_j$.
Since $[S_0,\Bbbk[p]]=0$,
the relation (\ref{3.3.4}) implies
$f_j(p)u_j=0$  for every $j=1,\dots,N$.

In particular, the matrix $P=(u_1, \dots, u_N)$ has non-zero
determinant.
In general, $P^{-1}$ does not lie in $M_N(\Bbbk[p])$,
but we still have $U_0=P\Bbbk^N$, $P^{-1}S_0 P = M_N(\Bbbk)$.
Then (\ref{3.3.3}) allows to conclude that
$P^{-1}M_N(\Bbbk[p]) P = M_N(\Bbbk[p])$.
The last relation implies
$P=g(p)R$, where $R,R^{-1}\in M_N(\Bbbk[p])$,
$g(p)\in \Bbbk [p]$.
For the matrix $R$ obtained we also have
$R^{-1}S_0 R= M_N(\Bbbk)$.

Let $\tau $ be the automorphism of $M_N(\Bbbk[p])$ defined
as the conjugation by~$R$. Then the composition
$\theta_1 = \tau\circ \theta^{-1}$ preserves $p$ and
$M_N(\Bbbk)$. Hence, $\theta_1 $ acts as the conjugation
by a nondegenerate matrix, i.e., $\theta_1: x\mapsto TxT^{-1}$
for some $T\in M_N(\Bbbk)$. Finally,
$\theta $ is the conjugation by invertible matrix
$Q= RT \in M_N(\Bbbk[p])$.

(ii)
Since $\partial_q(a)=[a,p]$ for every $a\in M_N(W)$,
we may conclude that $\theta $ satisfies the conditions
of Lemma~\ref{lem 3.3.3}.

It follows from (i) that
there exists an invertible matrix $Q\in M_N(\Bbbk[p])$ such that
the restriction $\theta\vert_{M_N(\Bbbk[p])}$
is just the conjugation by~$Q$.

If $S_0 = \theta^{-1}(M_N(\Bbbk))$, then
for $q+A = \theta^{-1}(q)$ we have
$[q+A,S_0]=0$. So $[Q^{-1}(q+A)Q, M_N(\Bbbk)]=0$,
and the matrix
$Q^{-1}(q+A)Q = q + Q^{-1}A Q + Q^{-1}Q'$ has to be scalar.
In particular,
$Q^{-1}AQ + Q^{-1}Q' = h(p)$ for some polynomial $h(p)$.
Hence,
$A = h(p) - Q'Q^{-1}$.
But $\theta(A) = Q^{-1}AQ$ by (i),
so $\theta(q)= q - \theta(A) = q- (h -Q^{-1}Q')$.
\end{proof}

\begin{thm}\label{thm 3.3.5}
Let $\theta $ be a $\partial_q$-invariant automorphism
 of $M_N(W)$. Then
\[
\theta = \theta_{\alpha, Q, h} : a(p,q) \mapsto
  Q^{-1}(p) a(p+\alpha, q - h(p)) Q(p)
\]
for some $\alpha \in \Bbbk$, $h(p)\in \Bbbk[p]$,
$Q,Q^{-1}\in M_N(\Bbbk[p])$.

The automorphism $\theta_{\alpha, Q, h}$ is continuous
if and only if $h(p)=0$.
\end{thm}

\begin{proof}
By Lemma \ref{lem 3.3.3} there exists $\alpha \in \Bbbk $
such that $\theta(p)=p+\alpha $.
Consider
$\tau_1 = \theta_{\alpha, \Id_N, 0}$, then
the composition
$\theta_1 = \tau_1^{-1}\circ \theta $ preserves~$p$:
$\theta_1(p)=p$.

Lemma \ref{lem 3.3.4} implies that there exist
a polynomial $h\in \Bbbk[p]$ and an invertible matrix
$Q\in M_N(\Bbbk[p])$
such that
$\theta_1 (q) = q - h(p) + Q^{-1}Q' $
and
$\theta_1\vert_{M_N(\Bbbk[p])}$
is the conjugation by~$Q$.
Consider $\tau_2 = \theta_{0,Q,h}$, then
$\tau_2 (q) = \theta_1 (q)$.
Since $\tau_2$ preserves $p$, we may conclude that
the composition $\theta_2 = \tau_2^{-1}\circ \theta_1$
is just the identity map.
So, $\theta = \tau_1\circ \tau_2 = \theta_{\alpha, Q, h}$.

It is clear that
$\theta_{\alpha, Q, 0} $
is continuous in the sense of $q$-topology.
If $h\ne 0$, then $\theta_{\alpha, Q, h}$
is not continuous since the sequence
$\theta_{\alpha, Q, h} (q^n)=Q^{-1}(q-h(p))^n Q $
does not converge to zero.
\end{proof}

\begin{cor}[{\cite[Theorem~4.1]{8}}]\label{cor 3.3.6}
Any automorphism of $\Cend_N$ is of the form $\Theta_{\alpha,Q}$,
as in (\ref{3.3.2}).
\end{cor}

\begin{proof}
If $\Theta $ is an automorphism of $\Cend_N$, then
the corresponding automorphism $\theta $ is of the form
$\theta_{\alpha, Q, 0}$ by
Proposition \ref{prop 3.3.1}
and Theorem \ref{thm 3.3.5}.
By the construction,
\[
\Theta(a)(n) = \theta_{\alpha, Q, 0}(a(n)) =
\Theta_{\alpha, Q}(a)(n) ,
\]
so $\Theta = \Theta_{\alpha, Q}$.
\end{proof}

\begin{prop}\label{prop 3.3.7}
Let $C \subseteq \Cend_N$ be a conformal subalgebra,
and let
$\theta: S(C) \to S$ be an isomorphism from
$S(C)$ onto a subalgebra $S\subseteq M_N(W)$.
If $\theta $ is $\partial_q$-invariant and continuous,
then there exists an injective homomorphism
$\Theta  : C \to {\mathcal R}(S)$
of
conformal subalgebras.
\end{prop}

\begin{proof}
Since $\theta $ is $\partial_q$-invariant and continuous,
the subalgebra $S$ satisfies the conditions
of Proposition \ref{prop 3.2.4}.

For every element $a\in C$ the sequence
$\{\theta(a(n))\}_{n=0}^\infty $ is differential.
By Lemma \ref{lem 3.2.3}
there exists an element $\hat a \in \Cend_N$
(moreover, $\hat a\in {\mathcal R}(S)$)
such that $\hat a(n) = \theta(a(n))$ for all $n\ge 0$.
The map $\Theta : a\mapsto \hat a$ is an injective
homomorphism of conformal algebras.
\end{proof}


\subsection{Left and right ideals of $\Cend_N$}\label{subsec3.4}

In this section, we state some results of \cite[Section~1]{8}
on the structure of left and right ideals of $\Cend_N$.

The following description of one-sided ideals was obtained
by B.~Bakalov.

\begin{prop}[\cite{8}]\label{prop 3.4.1}
The conformal algebra $\Cend_N$ is simple.
Every right ideal of $\Cend_N$ has
the form $\Cend_{P,N}=P(v)\Cend_N$,
where $P\in M_N(\Bbbk[v])$.
Every left ideal has the form
$\Cend_{N,Q}=\Cend_N \varphi^{-1}(Q)=
 \Cend_N Q(-D\otimes 1+1\otimes v)$,
where $Q \in M_N(\Bbbk[v])$, $\varphi$ is the isomorphism from
Proposition~\ref{prop 2.3.2}.
\quad \qed
\end{prop}

\begin{lem}[\cite{8}]\label{lem 3.4.2}
Every $\Cend_{N,Q}$,
$\det Q\ne 0$, is a conjugate of
$\Cend_{N,D}$, where
$D$ is the canonical diagonal form of $Q$, i.e.,
$D=\diag(f_1,\dots, f_N)$,
every $f_i$ is monic, and  $f_i \mid f_{i+1}$.
\quad \qed
\end{lem}

\begin{prop}[\cite{8}]\label{prop 3.4.3}
Conformal subalgebra
$\Cend_{N,Q}\subseteq \Cend_N$
is irreducible
if and only if $\det Q\ne 0$.
\quad \qed
\end{prop}

\begin{exmp}\label{exmp 3.4.4}
Consider $C = M_N(\Bbbk[D])\subset \Cend_N$.
This subalgebra
is just the current conformal algebra
$\Curr_N =\Curr M_N(\Bbbk)$.
It is clear (see also \cite{8})
that $\Curr _N$ is irreducible
since
$S(\Curr_N)=M_N(\Bbbk[q])$
and
$\Bbbk[p]S(\Curr_N)=M_N(W)$.
\end{exmp}

It was conjectured in \cite{21} that
$\Theta_{0,T}(\Curr_N)$, $T,T^{-1}\in  M_N(\Bbbk[v]) $,
and
$\Cend_{N,Q}$, $Q\in M_N(\Bbbk[v])$, $\det Q\ne 0$,
exhaust all irreducible subalgebras of $\Cend_N$.
In Section~\ref{sec5}, we will prove the conjecture.

\begin{prop}[\cite{8}]\label{prop 3.4.5}
If $\det Q\ne 0$, then:

{\rm (i)} every left ideal of $\Cend_{N,Q}$ is of the form
 $\Cend_{N,TQ}$, $T\in  M_N(\Bbbk[v])$;

{\rm (ii)} $\Cend_{N,Q}$ is a simple conformal algebra.
\quad \qed
\end{prop}

\begin{rem}\label{rem 3.4.6}
The structure of right ideals of
$\Cend_{N,Q}$, $\det Q\ne 0$,
was also described in \cite{8}.
In fact, every right ideal of $\Cend_{N,Q}$
is of the form $P\Cend_{N,Q}$, $P\in  M_N(\Bbbk[v])$.
\end{rem}


\section{Properties of operator algebras}\label{sec4}

\subsection{Density of irreducible subalgebras}\label{subsec4.1}

Proposition \ref{prop 3.2.6} shows the relation between irreducible
conformal subalgebras of $\Cend_N$ and irreducible
subalgebras of $M_N(W)$: if $C\subseteq \Cend_N$ is
irreducible, then $S_1=\Bbbk[p]S(C)\subseteq M_N(W)$
acts on $V_N$ irreducibly. The examples of irreducible
conformal subalgebras provided by
Proposition \ref{prop 3.4.3}
  and Example \ref{exmp 3.4.4}
correspond
to dense (in the sense of $q$-topology) subalgebras
of $M_N(W)$. In this section,
we prove that the density of $\Bbbk[p]S(C)$
is a necessary property of an irreducible conformal
subalgebra~$C$.

Let us remind the classical Density Theorem by N.~Jacobson.
Consider a primitive algebra $S$ over a field $\Bbbk$ with a faithful
irreducible left $S$-module~$V$. One may identify
$S$ with its image in $\End_{\Bbbk} V$.
The centralizer
 ${\mathcal D} = \{\phi \in \End_{\Bbbk} V\mid [\phi ,S]=0\}$
is a skew field over $\Bbbk $. The module $V$ could
be considered as a right vector space over
${\mathcal D}^{\mathrm{op}}$,
where
${\mathcal D}^{\mathrm{op}}$
is anti-isomorphic to~${\mathcal D}$.
There exists the natural embedding of $S$
into ${\mathcal E} = \End_{{\mathcal D}^{\mathrm{op}}} V$.

\begin{thm}[\cite{19}]\label{thm 4.1.1}
The algebra $S$ is a dense subalgebra of
${\mathcal E}$ in the sense of finite topology, i.e.,
for any ${\mathcal D}^{\mathrm{op}}$-independent family
$u_1,\dots, u_n\in V$ and for any $w_1,\dots, w_n\in V$
there exists $a\in S$ such that $au_i=w_i$, $i=1,\dots, n$.
\quad \qed
\end{thm}

\begin{thm}\label{thm 4.1.2}
Let $S\subseteq M_N(W)$ be a $\partial_q$-invariant
subalgebra such that $pS\subseteq S$.
If $S$ acts on $V_N$ irreducibly, then $S$ is a dense subalgebra
of $M_N(W)$ in the sense of $q$-topology.
\end{thm}

\begin{proof}
Let us fix $0\ne a\in S$.
By ${\mathcal D}\subseteq \End_{\Bbbk} V_N$
we denote the centralizer of~$S$.
For any
$\phi \in {\mathcal D}$
the relations
$[\phi,a]=[\phi, pa]=0$
imply
\begin{equation}\label{4.1.1}
[\phi,p]a=0.
\end{equation}
But for every $b\in S$ we have
\[
[[\phi,p],b]=[[\phi,b],p] + [\phi,[p,b]]=0,
\]
since $[p,b]= - \partial_q b\in S$.
Hence,
$[\phi, p]\in {\mathcal D}$
is either invertible or zero.
The first is not the case because of (\ref{4.1.1}), so
$[\phi, p]=0$.

Let us fix an $H$-basis $\{e_1,\dots, e_N\}$
of $V_N$.
Since $[\phi,p]=0$,
the matrix of $\phi $
in the $\Bbbk$-basis
$\{e_1,\dots, e_N, pe_1,\dots, pe_N, p^2e_1,\dots ,
 p^2e_N,\dots \}$
of $V_N$ has the block-triangular form
\[
[\phi ]=
\begin{pmatrix}
A_0 & 0 & 0 & \dots \\
* & A_0 & 0 & \dots \\
* & * & A_0 & \dots \\
\dots&\dots&\dots&\dots
\end{pmatrix},
\]
where
$A_0$
is a matrix from $M_N(\Bbbk)$.
Let $\lambda\in \Bbbk $ be an eigenvalue of $A_0$.
Then the transformation
$\phi - \lambda\Id_N$ lies  in ${\mathcal D}$,
and it is not invertible.
Hence,
$\phi=\lambda\Id_N$, $\lambda \in \Bbbk$.

Thus, ${\mathcal D} = \Bbbk $.
Theorem \ref{thm 4.1.1} implies $S$ to be a dense subalgebra
of $\End_{\Bbbk} V_N$ in the sense of finite topology.
It is clear that
the finite
topology is equivalent to the $q$-topology on $M_N(W)$,
so $S$ is dense in the sense
of $q$-topology.
\end{proof}


\subsection{Operator algebras of linear growth}\label{subsec4.2}

In this section, we prepare some additional
facts about operator algebras
of Gel'fand--Kirillov dimension one.
These subalgebras correspond to conformal
subalgebras of finite type.
We will essentially use the following result of \cite{27}.

\begin{thm}[\cite{27}]\label{thm 4.2.1}
Let $S$ be a finitely generated prime unital algebra
such that
$\GKdim S=1$. Then $S$ is a finitely generated module
over its center.
\quad \qed
\end{thm}

\begin{prop}\label{prop 4.2.2}
Assume that $S\subseteq M_N(W)$ is a $\partial_q$-invariant
subalgebra such that $S_1=\Bbbk[p]S$
is prime and satisfies ascending chain condition (a.c.c.)
for right annihilation ideals. Then $S$ is prime.
\end{prop}

\begin{proof}
Suppose that $S$ is not prime, i.e., there exist non-zero
ideals $I,J\ideal S$ such that $IJ=0$.

If either $I$ or $J$
is $\partial_q$-invariant, then
$(\Bbbk[p]I)(\Bbbk[p]J) = 0$,
where
$\Bbbk[p]I, \Bbbk[p]J \ideal_l S_1$.
Since $S_1$ is prime, either $I$ or $J$ is zero.

Consider the following ascending sequence of ideals of $S$:
$J_0=0$, $J_1=J$, $J_{n+1} = J_n+\partial_q(J_n)$ for $n\ge 1$.
It is easy to show by induction on~$n$ that $IJ_n \subseteq J_{n-1}$.
In particular, $I^n J_n = 0$ for every $n\ge 0$.

Since $S_1$ satisfies a.c.c. for right annihilation ideals,
there exists $m\ge 1$ such that
\[
\Ann_{S_1}(I)\subseteq \Ann_{S_1}(I^2)\subseteq
\dots
\subseteq \Ann_{S_1}(I^m) = \Ann_{S_1}(I^{m+1}) =\dots.
\]
We have shown that $J_n\subseteq \Ann_{S_1}(I^n)$, so
for every $p\ge m$
the ideal $J_p$ lies in $\Ann_{S_1}(I^m)$.
In particular, $I^m J_{\partial} = 0$, where
$J_{\partial} = \bigcup\limits_{n\ge 0} J_n$
is a $\partial_q$-invariant ideal.
If $I^m\ne 0$, then we
obtain a contradiction as it was shown above.

Therefore, it is sufficient to prove that $S$ has no
non-zero nilpotent ideals. By the very same reasons
as stated above, $S$ has no non-zero $\partial_q$-invariant
nilpotent ideals.

It is easy to note that, if $I^n=0$,
then $(I+\partial_q(I))^{n^2}=0$.
Hence, for every
nilpotent ideal $I\ideal S$ and for every $m\ge 0$
the ideal $I+\partial_q(I)+\dots + \partial_q^m(I)$
is also nilpotent. If $I\ne 0$, then we may choose
$0\ne a\in I$ and consider the ideal $I_a$
generated in $S$ by all derivatives of~$a$:
$I_a = (a, \partial_q(a),\dots ,\partial_q^m(a))\ideal S$.
This ideal is $\partial_q$-invariant and lies in
$I+\partial_q(I)+\dots + \partial_q^m(I)$, so $I_a$ is nilpotent.
As it was shown above, $I_a=0$ in contradiction with $a\ne 0$.
\end{proof}

For a subalgebra $S\subseteq M_N(W)$ denote
by $Z(S)$ the center of $S$.
If $S_1$ is a subalgebra of $M_N(W)$ such that
$S_1 \supseteq S$, then by
${\mathcal Z}_{S_1}(S)$
we denote the centralizer
of $S$ in $S_1$.
If $S$ and $S_1$ are $\partial_q$-invariant, then
so are ${\mathcal Z}_{S_1}(S)$ and $Z(S)$.

\begin{lem}\label{lem 4.2.3}
Let $S\subseteq M_N(W)$ be a finitely generated
prime
$\partial_q$-invariant
subalgebra such that
$\GKdim S=1$
and
${\mathcal Z}_{M_N(W)}(\Bbbk[p]S)=\Bbbk$.
Then there exists a matrix
$A\in M_N(\Bbbk[p])$ such that
$q+A\in Z(S)$.
\end{lem}

\begin{proof}
Although $S$ is not supposed to be unital, it is still
possible to apply Theorem \ref{thm 4.2.1}, since $S+\Bbbk$ also
satisfies the conditions of the statement.
By Theorem \ref{thm 4.2.1} $S$ is a finitely
generated module over its center
$Z(S)$. In particular, $Z(S)\ne 0,\Bbbk $.
Since $Z(S)$ is
$\partial_q$-invariant, there exists a non-zero element
$y \in Z(S)\cap M_N(\Bbbk[p]) $.
But $[y,p]=0$, so $y\in \Bbbk $ and $\Bbbk \subset Z(S)$. Now, let
$x\in Z(S)\setminus \Bbbk $. We may assume that
$\partial_q(x)\in \Bbbk$, hence,
$x=q+A$, where $A\in M_N(\Bbbk[p])$.
\end{proof}

\begin{prop}\label{prop 4.2.4}
Let $S\subset M_N(W)$ be a $\partial_q$-invariant
subalgebra such that
\begin{equation}\label{4.2.1}
\bigoplus\limits_{n\ge 0}p^nS = M_N(W).
\end{equation}
Then there exists a $\partial_q$-invariant
automorphism $\theta $ of $M_N(W)$ such
that $\theta(p)=p$ and $\theta(S) = M_N(\Bbbk[q])$.
\end{prop}

\begin{proof}
It is clear that $S$ is finitely generated and
$\GKdim S=1$ (see, e.g., \cite{26}).
Since $M_N(W)$ is (left and right) Noetherian \cite{15},
it satisfies a.c.c. for
right annihilation ideals. By Proposition \ref{prop 4.2.2}
the algebra $S$ is prime,
and by Lemma \ref{lem 4.2.3} there exists $x= q+A\in Z(S)$.

Moreover, for every $\partial_q$-invariant ideal $I\ideal S$
its ``envelope'' $\Bbbk[p]I$ is an ideal of $M_N(W)$. Thus,
either $I=0$ or $\Bbbk[p]I=M_N(W)$. The last relation implies
$I=S$ since the sum (\ref{4.2.1}) is direct. Therefore, $S$
has no non-trivial $\partial_q$-invariant ideals.

Denote $S_0 =S\cap M_N(\Bbbk[p])=\Kerr \partial_q\vert_S$.
For every $k\ge 1$ let us define
$S_k =\{a\in S \mid \partial_q(a)\in S_{k-1} \}$.
Since $\partial_q$ is locally nilpotent, we have
\[
S=\bigcup\limits_{k\ge 0} S_{k},
\quad S_{k}=\Kerr \partial_q^{k+1}\vert_S.
\]
It is easy to show by induction on $k$
that
\begin{equation}\label{4.2.2}
S_{k}=S_{0}+S_{0}x + \dots + S_{0} x^k.
\end{equation}
Indeed, it is clear for $k=0$. For every $k\ge 1$
and for every $a\in S_{k}$
we have
\[
\partial_q^k \left(a - \frac{1}{k}\partial_q(a)x \right) =0,
\]
so
$a - \frac{1}{k}\partial_q(a)x\in S_{k-1}$.
Since $\partial_q(a)\in S_{k-1}$, we may use the
inductive assumption. Hence, (\ref{4.2.2}) is true.

It is also clear that for
every $k\ge 0$ the sum (\ref{4.2.2}) is direct. So,
$S=\bigoplus\limits_{n\ge 0}S_{0}x^n$
and $\GKdim S_{0}=0$ \cite{26}.
The algebra $S_{0}$ has no non-trivial ideals since $S$ is
differentially simple. Moreover, for this simple
finite-dimensional algebra $S_0$ we have
\[
\Bbbk[p]S_0 = \bigoplus\limits_{n\ge 0}p^nS_0 = M_N(\Bbbk[p]).
\]
Hence,  there exists an isomorphism $\theta: S_0 \to M_N(\Bbbk)$
which could be
extended via $\theta(x)=q$, $\theta(p)=p$
to a $\partial_q$-invariant automorphism of $M_N(W)$
such that
$\theta(S) = M_N(\Bbbk[q])$.
\end{proof}


\subsection{Operator algebras of current-type
 subalgebras}\label{subsec4.3}

The main result of this section could be stated as follows:
if there exists a $\partial_q$-invariant automorphism of
$M_N(W)$ such that an operator algebra
$S=S(C)$ goes onto $M_N(\Bbbk[q])$,
then $C$ is a conjugate of $\Curr_N$
by an automorphism of $\Cend_N$
(i.e., $C$ is a {\em current-type\/} conformal subalgebra
of $\Cend_N$). To prove this, we will show that one
may choose the isomorphism
$\theta $
in Proposition \ref{prop 4.2.4} to be continuous.

Let $h\in \Bbbk[p]$ be a polynomial.
Consider two sequences of polynomials in~$p$:
$\{h_{(n)}\}_{n\ge 0}$ and
$\{h_{[n]}\}_{n\ge 0}$
such that
\begin{equation}\label{4.3.1}
h_{(0)}=h_{[0]}=1,\quad
h_{(n)}=-hh_{(n-1)}+ h_{(n-1)}',\quad
h_{[n]}= hh_{[n-1]}+ h_{[n-1]}'.
\end{equation}
It is easy to see that these sequences are constructed
by the following reason. For any $n\ge 0$
we may write
\begin{equation}\label{4.3.2}
(q-h)^n = a_n(p,q)q + h_{(n)},
\quad (q+h)^n=b_n(p,q)q + h_{[n]}
\end{equation}
in the Weyl algebra~$W$.

\begin{lem}\label{lem 4.3.1}
For every
$k\ge 0$ we have
\begin{equation}\label{4.3.3}
\sum\limits_{s=\xi}^k
 \binom{k-\xi}{s-\xi} h_{[s-\xi]}h_{(k-s)}=
 \begin{cases}
  0, & \xi<k, \\
        1,& \xi=k.
      \end{cases}
\end{equation}
\end{lem}

\begin{proof}
The case $\xi =k$ is obvious, so we will concentrate
our attention
on the case  $\xi<k$.
For $k=1$, $\xi= 0$ the relation (\ref{4.3.3}) is clear.
Suppose that (\ref{4.3.3}) holds for some integer
$k$ and for all $\xi\le k$.
Then for $k+1$ instead of $k$ and for $\xi =k$
the corresponding equality in (\ref{4.3.3}) is also true.
For $\xi <k$ we have
\begin{eqnarray}
&& \sum\limits_{s=\xi}^{k+1}
  \binom{k-\xi +1}{s-\xi} h_{[s-\xi]}h_{(k+1-s)}
  \nonumber
 \\
&&\qquad
   =
  \sum\limits_{s=\xi}^k \binom{k-\xi}{s-\xi} h_{[s-\xi]}
  \big( -hh_{(k-s)}+ h_{(k-s)}' \big)
   \nonumber
    \\
&& \qquad
  +
  \sum\limits_{s=\xi+1}^{k+1}\binom{k-\xi}{s-\xi -1}
  \big( hh_{[(s-1)-\xi]} + h_{[(s-1)-\xi]}' \big) h_{(k-(s-1))}
    \nonumber
   \\
&& \qquad
    =
   \sum\limits_{s=\xi}^k \binom{k-\xi}{s-\xi}
   \big( h_{[s-\xi]}h'_{(k-s)} + h'_{[s-\xi]}h_{(k-s)} \big)
   =0.
     \nonumber
\end{eqnarray}
\end{proof}

\begin{cor}\label{cor 4.3.2}
For any $\xi<k$ we have
\begin{equation}\label{4.3.4}
\sum\limits_{s=\xi}^{k-1} \binom{k}{s}\binom{s}{\xi}
 h_{[s-\xi]}h_{(k-s)} = -\binom{k}{\xi} h_{[k-\xi]}.
 \quad \qed
\end{equation}
\end{cor}

\begin{lem}\label{lem 4.3.3}
Let
\[
a(n)=\sum\limits_{k\ge 0}\binom{n}{k} A_k(p)q^{n-k}
\]
be an element of $M_N(W)$.
Replace $q$ with $x-h$, where $x=q+h$, and consider the presentation
\[
a(n)=\sum\limits_{s\ge 0}\binom{n}{s} B_s(p)x^{n-s}.
\]
Then
\begin{equation}\label{4.3.5}
A_k = \sum\limits_{s\ge 0}\binom{k}{s} B_s(p)h_{[k-s]}.
\end{equation}
\end{lem}

\begin{proof}
Since $\partial_x = \partial_q=[\cdot, p]$ on $M_N(W)$,
we may assume $A_0=B_0$,
$A_1 = B_0 h + B_1$. For greater $k$, the relation (\ref{4.3.5})
could be proved by induction.
Indeed,
\[
B_k=A_k + \sum\limits_{s=0}^{k-1} \binom{k}{s}
 A_s h_{(k-s)},
\]
so by the assumption
\[
A_k = B_k
  - \sum\limits_{s=0}^{k-1}\sum\limits_{\xi=0}^s
  \binom{k}{s}\binom{s}{\xi} B_\xi h_{[s-\xi]} h_{(k-s)},
\]
and it is sufficient to apply (\ref{4.3.4}).
\end{proof}

\begin{lem}\label{lem 4.3.4}
Let $C\subset \Cend_N$ be a conformal
subalgebra such that $S(C)=M_N(\Bbbk[q+h(p)])$.
Then $\deg h=0$ (or $h=0$).
\end{lem}

\begin{proof}
Consider an arbitrary element $a\in C$:
\[
a=\sum\limits_{k=0}^n (-D)^{(k)}\otimes A_k(v)\ne 0.
\]
Then
\begin{equation}\label{4.3.6}
a(n+1)=\sum\limits_{k=0}^n
 \binom{n+1}{k} A_k q^{n+1-k}\in S(C)
 =M_N(\Bbbk)[q+h(p)].
\end{equation}
We may assume that $A_0\ne 0$:
otherwise, $a=D^{(m-1)}b$, $b\in \Cend_N$, $m\ge 2$
and it is sufficient
to consider $a(n+m)$ instead of $a(n+1)$
(although the element $b$ may not lie in $C$, the operator $b(n+1)$
lies in $S(C)=M_N(\Bbbk [q+h(p)])$).
Rewriting (\ref{4.3.6}) by using $q=x-h$ gives
\[
B_{n+1}=\sum\limits_{k=0}^n \binom{n+1}{k} A_k h_{(n+1-k)}
 \in M_N(\Bbbk).
\]
But Lemma \ref{lem 4.3.3} and Lemma \ref{lem 4.3.1} imply that
\begin{eqnarray}
B_{n+1} & = & \sum\limits_{k,s=0}^n
  \binom{n+1}{k}\binom{k}{s} B_s h_{[k-s]}h_{(n+1-k)}
 \nonumber \\
  & = & \sum\limits_{s=0}^n \binom{n+1}{s}B_s
    \sum\limits_{k=s}^n \binom{n+1-s}{k-s} h_{[k-s]}h_{(n+1-k)}
    \nonumber \\
  & = &
    -\sum\limits_{s=0}^n \binom{n+1}{s} B_s h_{[n+1-s]}.
    \label{4.3.7}
\end{eqnarray}
If $\deg h>0$, then the term of highest degree in (\ref{4.3.7})
corresponds to $s=0$:
\[
B_{n+1} = -A_0 h_{[n+1]} + \hbox{terms of lower degree}
  \notin M_N(\Bbbk).
\]
Hence, either $h=0$ or $\deg h=0$.
\end{proof}

\begin{thm}\label{thm 4.3.5}
Let $S\subset M_N(W)$ be a subalgebra satisfying
 the conditions of Proposition \ref{prop 4.2.4}
and Proposition \ref{prop 3.2.4},
 i.e., $S=S(C)$ for a conformal
subalgebra $C\subset \Cend_N$. Then
there exists an automorphism $\Theta $ of the conformal
algebra $\Cend_N$ such that $\Theta(C) = \Curr_N$.
\end{thm}

\begin{proof}
By Proposition \ref{prop 4.2.4}, there exists a $\partial_q$-invariant
automorphism $\theta $ of $M_N(W)$ such that
$\theta(S)=M_N(\Bbbk[q])$, $\theta(q+A(p)) = q$.
The matrix $A(p)$ is defined up to an additive scalar, i.e.,
instead of $q+A(p)$ one may consider $q+A(p)+\alpha\in Z(S)$,
for any $\alpha \in \Bbbk$.
Theorem \ref{thm 3.3.5} implies that $\theta = \theta_{0,Q,h}$,
$A=h(p) - Q'_pQ^{-1}$, and we may assume that $h(0)=0$.

Consider the automorphism $\Theta_{0,Q}$ of  the conformal
algebra $\Cend_N$, and denote $\Theta_{0,Q}(C)=C^Q$. Then
the subalgebra
\[
S^Q = S(C^Q) = Q^{-1}S(C) Q \subset M_N(W)
\]
is also isomorphic to $M_N(\Bbbk[q])$,
and the isomorphism
is given by $\theta_{0,\Id_N,h}$, where $h(0)=0$.

By Lemma \ref{lem 4.3.4}  $h=0$, so $\theta =\theta_{0,Q,0}$
is continuous by Theorem \ref{thm 3.3.5}.
Proposition \ref{prop 3.3.7}
and Corollary \ref{cor 3.2.9}
imply  $C^Q\ideal_l {\mathcal R}(M_N(\Bbbk[q]))= \Curr_N$.
It is well-known
 that every left
ideal of $\Curr_N$ is of the form $H\otimes I_0$,
where $I_0\ideal_l M_N(\Bbbk)$.
If $C^Q = H\otimes I_0$, then $S^Q=M_N(\Bbbk[q])=I_0\otimes \Bbbk[q]$;
hence, $I_0=M_N(\Bbbk)$.
Thus, we obtain $C^Q = \Curr_N$.
\end{proof}


\section{Irreducible conformal subalgebras}\label{sec5}

\subsection{Preliminary notes}\label{subsec5.1}

\begin{thm}\label{thm 5.1.1}
Let $C\subseteq \Cend_N$ be an irreducible conformal subalgebra.
Then
\begin{equation}\label{5.1.1}
\Bbbk[v]C =
C+vC+v^2C+\dots = \Cend_{N,Q},
\end{equation}
where $Q=Q(v)\in M_N(\Bbbk[v])$ is a nondegenerate matrix.
\end{thm}

\begin{proof}
By Proposition \ref{prop 3.2.6}, the subalgebra
$S_1=\Bbbk[p]S(C)$ acts irreducibly on $V_N$.
So Theorem \ref{thm 4.1.2}
implies $S_1$ to be a dense subalgebra of $M_N(W)$
in the sense of $q$-topology.
In particular, for every
$A(p)q^n\in M_N(\Bbbk[p])q^n$,
and for every integer $M\ge n+1$
there exists
$a_{A,n, M}\in S_1$
such that
$a_{A,n, M} - A(p)q^n \in M_N(W) q^M $,
i.e.,
\[
a_{A,n, M} = A(p)q^n + x(p,q)q^M
\]
for some $x(p,q)\in M_N(W)$.

The $H$-submodule $C_1=\Bbbk[v]C \subseteq \Cend_N$
is also a conformal subalgebra of $\Cend_N$
(see (\ref{3.1.1})), and $S(C_1)=S_1$.
So Corollary \ref{cor 3.2.8}
implies $S_1 \cdot C_1\subseteq C_1$.

Let us fix an element $b\in C_1$:
\[
b=\sum\limits_{s= 0}^m D^{(s)}\otimes B_s.
\]
For any $n\ge 0$, choose an integer
 $M > n+ m + \max\limits_{s=0,\dots,m}\deg_v(B_s)$.
Then by Proposition~\ref{prop 3.2.7}
\[
a_{A,n,M}\cdot b = A(v)\oo{n} b\in C_1
\]
for every $A(p) q^n\in M_N(W)$.
Hence, $C_1$ is a left ideal of $\Cend_N$, and by
Proposition \ref{prop 3.4.1}
it is of the form
$\Cend_{N,Q}$
for some $Q\in M_N(\Bbbk[v])$.
Proposition \ref{prop 3.4.3}
implies $\det Q \ne 0$.
\end{proof}

\begin{lem}\label{lem 5.1.2}
{\rm (i)}
If $I$ is a left (right) ideal of $C$, then
$\Bbbk[v]I$ is a left (right) ideal of $\Bbbk[v]C=\Cend_{N,Q}$.

{\rm (ii)}
If $0\ne I\ideal C$, then $\Bbbk[v]I=\Cend_{N,Q}$.
\end{lem}

\begin{proof}
Statement (i) easily follows from (\ref{3.1.1}).

To prove (ii), it is sufficient to note that
$\Bbbk[v]I\ne 0$ is an ideal of $\Cend_{N,Q}$
for every non-zero ideal $I$ of $C$.
By Proposition \ref{prop 3.4.5}(ii),
the conformal algebra
$\Cend_{N,Q}$ is simple, so $\Bbbk[v]I=\Cend_{N,Q}$.
\end{proof}

\begin{prop}\label{prop 5.1.3}
If there exists an element $a\in C$ such that
$a\ne 0$ and $v^ka\in C$ for all $k\ge 0$, then
$C=\Cend_{N,Q}$.
\end{prop}

\begin{proof}
It is clear from (\ref{3.1.1}) that
$I=\{a\in C\mid v^k a \in C,\, k\ge 0 \}$
is an ideal of $C$. In particular, $\Bbbk[v]I=I$.
If $I\ne 0$, then
by Lemma \ref{lem 5.1.2}(ii) we have
$C\supseteq I=\Bbbk[v]I=\Cend_{N,Q}$.
\end{proof}

For every irreducible subalgebra $C\subseteq \Cend_N$
there exist three options as follows:

Case 1: the sum (\ref{5.1.1}) is direct, i.e.,
\begin{equation}\label{5.1.2}
\Cend_{N,Q}=\bigoplus\limits_{n\ge 0}v^nC.
\end{equation}

Case 2: $C\cap vC\ne 0$.

Case 3: the sum (\ref{5.1.1}) is non-direct, but $C\cap vC=0$.

We will show that the first case corresponds to the
current-type conformal subalgebras,
the second one gives $C=\Cend_{N,Q}$, and the third one is impossible.

Without loss of generality we may assume $Q$ to be in the canonical
diagonal form: $Q=\diag(f_1,\dots, f_N)$, where $f_i$ are monic
polynomials and $f_i\mid f_{i+1}$.
Indeed, for an arbitrary
$Q\in M_N(\Bbbk[v])$, $\det Q\ne 0$, there exist
$U, T\in M_N(\Bbbk[v])$,
$\det U, \det T\in \Bbbk\setminus\{0\}$,
such that $TQU = D$, where $D$ has the canonical diagonal form.
Then for $C^U = \Theta_{0,U}(C)$ we have
$\Bbbk[v]C^U = \Cend_{N,D}$, and all the conditions described by
Cases 1--3 hold.


\subsection{Finite type case}\label{subsec5.2}

Now, let us consider Case~1.
Throughout this section $C$
is a conformal subalgebra of $\Cend_N$ satisfying
(\ref{5.1.2}) for a fixed matrix
$Q=\diag(f_1,\dots, f_N)$,
$0\ne f_i\in \Bbbk[v]$, $i=1,\dots, N$.

\begin{lem}\label{lem 5.2.1}
{\rm (i)}
If $h(D)a\in \Cend_{N,Q}$ for some $0\ne h\in H$,
then $a\in \Cend_{N,Q}$.

{\rm (ii)}
For $S=S(C)$
we have
\begin{equation}\label{5.2.1}
M_N(W)Q=\bigoplus\limits_{n\ge 0} p^n S.
\end{equation}
\end{lem}

\begin{proof}
(i) Consider $U = \Cend_N/\Cend_{N,Q}$ as a module over
the conformal algebra $\Cend_N$. If $h(D)a\in \Cend_{N,Q}$,
then $\bar a= a+\Cend_{N,Q}$ lies in the torsion of $U$.
In particular,
$\Cend_N \oo{\omega } a \in \Cend_{N,Q}$. But there exists
$e=1\otimes \Id_N\in \Cend_{N}$ such that
$e\oo{0} a = a \in \Cend_{N,Q}$.

(ii)
Suppose that the sum (\ref{5.2.1}) is non-direct.
Since $(Da)(n) = -na(n-1)$, we may consider
\[
0=a_0(n) + \dots + p^ma_m(n),\quad  a_i\in C,\ n\ge 0.
\]
Then for $a=a_0+va_1+\dots +v^ma_m\in \Cend_{N,Q}$
we obtain $a(0)=\dots=a(n)=0$, i.e., $a=D^{n+1}b$.
The statement (i) implies $b\in \Cend_{N,Q}$, so
$b=b_0+vb_1+ \dots +v^kb_k$, $b_i\in C$.
Since (\ref{5.1.1}) is direct, we have
$k=m$ and $D^{n+1}b_i=a_i$, so $a_i(n)=0$ for all $i=0,\dots, m$.
\end{proof}

Therefore, the operator algebra $S=S(C)$ of the initial
conformal subalgebra $C$ satisfies (\ref{5.2.1}).
Consider the presentations of all elements
\begin{eqnarray}
& q^m e_{ij}Q = \sum\limits_{s\ge 0} p^s a_{ijm, s},
 \nonumber \\
& i,j=1,\dots , N,\ m=0,\dots, \max\limits_{j=1,\dots, N}\deg f_j+1,
\ a_{ijm,s}\in S,
  \nonumber
\end{eqnarray}
where $e_{ij}$ are the matrix units
(so $e_{ij}Q=e_{ij}f_j$).
The finite set  of all $a_{ijm,s}$ together with all
their derivatives
$\partial_q^n a_{ijm,s}$ generates
a $\partial_q$-invariant subalgebra
$S_0$ of $S$.
Denote $S_{01}=\Bbbk[p]S_0$,
then $S_{01}$ is a $\partial_q$-invariant subalgebra of $S_1=M_N(W)Q$.

\begin{lem}\label{lem 5.2.2}
{\rm (i)} $\GKdim S_{01}=2$;

{\rm (ii)} ${\mathcal Z}_{M_N(W)}(S_{01})=\Bbbk$;

{\rm (iii)} $S_{01}$ has no non-zero nilpotent ideals.
\end{lem}

\begin{proof}
(i) Since $\GKdim M_N(W)=2$, we have
$\GKdim S_{01}\le 2$.
Let us consider the subalgebra $S_{00}$ generated by
$f_1e_{11}$, $qf_1e_{11}$ in $S_{01}$.
If $\deg f_1>0$, then $\GKdim S_{00}=2$, so
$\GKdim S_{01}=2$.
If $\deg f_1=0$, then we obtain that $We_{11}\subseteq S_{01}$,
hence, $\GKdim S_{01}=2$ as well.

(ii)
Let $a\in {\mathcal Z}={\mathcal Z}_{M_N(W)}(S_{01})$.
Then, in particular,
$[e_{ij}f_j, a]=0$, so $a$ is a diagonal matrix.
Moreover,
$[pe_{ij}f_j,a]=0 = [p,a]e_{ij}f_j$
for all $i,j=1,\dots, N$.
Hence, $[p,a]=0$.
By the very same reasons $[a,q]=0$.
Since all $f_j$ are unital, we conclude that $a=\alpha \Id_N$,
$\alpha\in \Bbbk$.

(iii)
If $0\ne I \ideal S_{01}$,
then for any $a\in I$ and for all $i,j,k,l=1,\dots, N$
the element $f_je_{ij} a f_k e_{lk}$ lies in $I$. In particular,
$I$ contains an element
$g e_{11}$
for some $0\ne g\in W$. Then
$I^n\ni g^ne_{11}\ne 0$ for all $n\ge 1$.
\end{proof}

\begin{lem}\label{lem 5.2.3}
{\rm (i)} $S_{01} \Bbbk[q] = M_N(W)$;

{\rm (ii)} if $I\ideal_r S_{01}$, then $I\Bbbk[q]\ideal_r M_N(W)$.
\end{lem}

\begin{proof}
(i)
If
$Q=\Id_N$ (i.e.,
$\deg f_1=\dots = \deg f_N =0$),
then the statement is clear.
If $\deg f_j>0$ for some $j\in \{1,\dots, N\}$,
then consider
\begin{eqnarray}
&  S_{01}\ni q^kf_je_{ij} = f_jq^ke_{ij} +
 k \partial_p(f_j)q^{k-1}e_{ij} +\dots
  + \partial_p^k(f_j) e_{ij},
   \nonumber \\
&  k=1,\dots, \deg f_j.
   \nonumber
\end{eqnarray}
It is clear that $\partial_p^k(f_j) e_{ij}\in S_{01}\Bbbk[q]$
for all $k=1,\dots, \deg f_j$,
so, $e_{ij}\in S_{01}\Bbbk[q]$.

(ii) It follows from (i).
\end{proof}

\begin{prop}\label{prop 5.2.4}
The algebra $S_{01}$ is prime and satisfies a.c.c. for
right annihilation ideals.
\end{prop}

\begin{proof}
First, let us show that $S_{01}$ satisfies a.c.c. for right
ideals of the form
$\Ann_{S_{01}}(X) = \{a\in S_{01}\mid Xa=0\} $,
$X\subseteq S_{01}$.

Consider  an ascending chain
of right annihilation ideals of $S_{01}$:
\[
I_1\subseteq I_2 \subseteq \dots , \quad
I_k = \Ann_{S_{01}}(X_k),\quad  X_k\subseteq S_{01}.
\]
By Lemma \ref{lem 5.2.3}(ii),
$I_1\Bbbk[q]\subseteq I_2\Bbbk[q]\subseteq \dots $
is an ascending chain of right ideals of $M_N(W)$.
Since $M_N(W)$ is a Noetherian algebra
 \cite{15}, we have
$I_n\Bbbk[q]=I_{n+1}\Bbbk[q]$ for a sufficiently large number~$n$.
In particular, $I_{n+1}\subseteq I_n\Bbbk[q]$.
Thus,
$X_{n}I_{n+1}\subseteq X_nI_n\Bbbk[q]= 0$, so
$I_{n+1}=I_n$.

Now, suppose that there exist two ideals
$I,J \ideal S_{01}$ such that $IJ=0$.
If $I\ne 0$, then by Lemma \ref{lem 5.2.2}(iii)
 the ideal $I$ is not nilpotent.
By the very same reason
as in the proof of Proposition~\ref{prop 4.2.2}, we may assume that
$J$ is $\partial_q$-invariant.

If $J$ is a non-zero $\partial_q$-invariant ideal of $S_{01}$,
then by Lemma \ref{lem 5.2.3}(ii)
 $J_1=J\Bbbk[q]\ideal_r M_N(W)$ is also
$\partial_q$-invariant.  Moreover, $S_{01}J_1\subseteq J_1$.
Then for every pair $(i,j)$, $i,j=1,\dots, N$,
there exist $k,l\in \{1,\dots, N\}$
such that
$e_{ik}f_k J_1 e_{lj} f_j \ne 0$.
Thus,
we have
$J_{ij}=J_1\cap e_{ij}W \ne 0$
for all $i,j=1,\dots, N$.
It is clear that
$A_{ij}=\{a\in W \mid ae_{ij}\in J_{ij}\}\ideal_r W$,
and $\partial_q(A_{ij})\subseteq A_{ij}$.
Then $A_{ij}=g_{ij}(p)W$
for some polynomial $0\ne g_{ij}\in \Bbbk[p]$.
In particular, $P=\sum\limits_{i=1}^N g_{ii}e_{ii}\in J_1$
is a diagonal matrix such that $\det P\ne 0$.

The relation $IJ=0$ implies $IJ_1 = (IJ)\Bbbk[q] = 0$,
so $IP=0$. Then $I=0$ since $\det P\ne 0$.
\end{proof}

\begin{cor}\label{cor 5.2.5}
$Q=\Id_N$,
$\Cend_{N,Q}=\Cend_N$,
\end{cor}

\begin{proof}
 By the construction, $S_{01}=\bigoplus\limits_{n\ge 0}p^n S_0$,
 so $\GKdim S_0=1$. Lemma \ref{lem 4.2.3} implies that there exists
 a matrix $A\in M_N(\Bbbk[p])$ such that $q+A\in Z(S_0)$.
 Then $\partial_p(x) = [x,A]$ for every $x\in S_0$.

Let us consider
\[
S_{01}\ni Q=a_0 + pa_1 +\dots + p^n a_n,\quad a_i\in S_0.
\]
Then
\[
\partial_p(Q)  = \partial_p(a_0) +p\partial_p(a_1) + \dots +
p^n\partial_p(a_n) + a_1 +\dots + np^{n-1}a_n.
\]
On the other hand,
\[
[Q,A(p)]= \partial_p(a_0) +p\partial_p(a_1) + \dots +
p^n\partial_p(a_n).
\]
So, $\partial_p(Q) - [Q,A(p)] \in S_{01}\subseteq M_N(W)Q$.
It is easy to note that it is not possible, if $\deg f_j\ne 0$
for at least one $j\in \{1,\dots, N\}$.
\end{proof}

\begin{thm}[{\cite[Theorem~5.2]{8}}]\label{thm 5.2.6}
There exists an automorphism $\Theta=\Theta_{0,P}$
of $\Cend_N$ such that
$\Theta(C) = \Curr_N $.
\end{thm}

\begin{proof}
By Lemma \ref{lem 5.2.1}(ii) and Corollary \ref{cor 5.2.5},
the operator algebra $S=S(C)$ satisfies the conditions
of Theorem \ref{thm 4.3.5}. The last statement implies
that $C$ is a conjugate of $\Curr_N$.
\end{proof}


\subsection{The case $C\cap vC\ne 0$}\label{subsec5.3}

Let $C$ be an irreducible conformal subalgebra of $\Cend_N$
such that $C\cap vC\ne 0$.

Remind the definition of the map $\varphi $
from Proposition~\ref{prop 2.3.2}.
For any element $a=f(D)\otimes A(v)\in H\otimes M_N(\Bbbk[v])$
we put
\begin{eqnarray*}
\varphi(a) & = & (f(D)\otimes 1)A(D\otimes 1 + 1\otimes v),  \\
\varphi^{-1}(a) & = & (f(D)\otimes 1)A(-D\otimes 1 + 1\otimes v). \\
\end{eqnarray*}

First, let us consider the main technical
features of \cite{8} concerning the
description of irreducible subalgebras for $N=1$.

\begin{lem}[c.f. {\cite[Section~2]{8}}]\label{x-lem5.3.1}
The following relations hold in $\Cend_N$:
\begin{eqnarray}
& (1\otimes A)\oo{n} \varphi^{-1}(B)
   =
 \begin{cases}
  (1\otimes A)\varphi^{-1}(B), & n=0, \\
           0, & n>0,
           \end{cases}
                                \label{eq-x5.3.1}\\
& (1\otimes A)\varphi^{-1}(B) \oo{n} (1\otimes A_1)
    =
    1\otimes A \partial_v^n (BA_1) ,
                           \label{eq-x5.3.3} \\
& (1\otimes A)\varphi^{-1}(B) \oo{n} (1\otimes A_1)\varphi^{-1}(B_1)
   =
  (1\otimes A\partial_v^n(BA_1)) \varphi^{-1}(B_1),
                             \label{eq-x5.3.4}
\end{eqnarray}
where $A,B,A_1, B_1 \in M_N(\Bbbk [v])$.
\end{lem}

\begin{proof}
Note that
\[
\varphi^{-1}(B) =
  \sum\limits_{s\ge 0} (-D)^{(s)}\otimes \partial_v^s(B),
\]
so (\ref{eq-x5.3.1}) follows from the axiom (C3).

Relation (\ref{eq-x5.3.3}) follows from (C2) and the Leibnitz
formula. In order to obtain (\ref{eq-x5.3.4}) one may use
(\ref{eq-x5.3.1}), (\ref{eq-x5.3.3}) and the associativity
relation~(\ref{2.1.3}).
\end{proof}

\begin{thm}[\cite{8}]\label{thm 5.3.1}
Let $C\subseteq \Cend_1$ be an irreducible
conformal subalgebra, and let $C\cap v\Cend_1 \ne 0$.
Then $C=\Cend_{1,Q}$ for some $0\ne Q\in \Bbbk [v]$.
\end{thm}

\begin{proof}
By Theorem \ref{thm 5.1.1}, $\Bbbk[v]C=\Cend_{1,Q}$.
Let $0\ne a\in C$ be an arbitrary element.
Let us present it as
$a=\sum\limits_{s=0}^m D^{s}\otimes a_s$. Consider
$\varphi(a)$, where $\varphi $ is the isomorphism from
Proposition \ref{prop 2.3.2}. Present
$b=\varphi(a)=\sum\limits_{t=0}^n D^{t}\otimes b_t$.
Here $a_m, b_n\in \Bbbk[v]\setminus \{0\}$.
Then
\[
a=\varphi^{-1}(b)=\sum\limits_{t=0}^n
(D^{t}\otimes 1) \varphi^{-1}(b_{t}).
\]
It follows from (\ref{eq-x5.3.1}), (C2), (C3) that
${\mathcal N}(a,a)=n+m+1$
and
\begin{eqnarray}
a\oo{n+m} a & = & a \oo{n+m} \varphi^{-1}(b)
                \nonumber \\
  & = &
   \left (
   \sum\limits_{s=0}^m D^{s}\otimes a_s
   \right)
     \oo{n+m}
   \left (
   \sum\limits_{t=0}^n
    (D^{t}\otimes 1)  \varphi^{-1}(b_{t})
   \right)
    \nonumber \\
  & = &
   (-1)^m (m+n)! (1\otimes a_m )  \varphi^{-1}(b_n).
     \nonumber
\end{eqnarray}
So, $(1\otimes a_m )\varphi^{-1}(b_n)\in C$. Since
$C\cap v\Cend_1 \ne 0$, we may assume
$\deg_v a_m>0$.

Therefore, $C$ contains an element
of the form
$(1\otimes g)  \varphi^{-1}(f)$,
$f,g\in \Bbbk [v]\setminus \{0\}$,
$\deg_v g>0$.
Consider two elements of this form in $C$
(they might be equal):
$a = (1\otimes g ) \varphi^{-1}(f)$,
$b = (1\otimes g ) \varphi^{-1}(h)$.
It follows from (\ref{eq-x5.3.4}) that
 ${\mathcal N}(a,b)= n+1$, where
$n= \deg_v f +\deg_v g $, and
\begin{eqnarray}
&  a\oo{n} b = \gamma  (1\otimes g )  \varphi^{-1}(h), \quad
 \gamma \in \Bbbk \setminus \{0\}, \nonumber \\
&  a\oo{n-1} b = \gamma'  (1\otimes (v+\alpha)g )
    \varphi^{-1}(h), \quad
 \gamma' \in \Bbbk \setminus \{0\},
 \
 \alpha \in \Bbbk.
 \nonumber
\end{eqnarray}
So we may conclude that $
(1\otimes vg ) \varphi^{-1}(h)  \in C$
as well as $(1\otimes v^k g )  \varphi^{-1}(h) \in C$ for all
$k\ge 0$.
Hence, the conditions of Proposition \ref{prop 5.1.3}
 hold
and $C=\Cend_{1,Q}$.
\end{proof}

We will use this technique in the general case.

\begin{thm}\label{thm 5.3.2}
Let $C\subseteq \Cend_N$ be an irreducible
conformal subalgebra, and let $C\cap vC\ne 0$.
Then $C=\Cend_{N,Q}$ for some $Q\in M_N(\Bbbk[v])$,
$\det Q\ne 0$.
\end{thm}

\begin{proof}
By Theorem \ref{thm 5.1.1}, $\Bbbk[v]C=\Cend_{N,Q}$,
$\det Q\ne 0$.
Let $I_1=\{a\in C \mid va\in C\}\ne 0$.
It is clear that $I_1$  is a non-zero ideal of~ $C$.
Moreover, $I_2=\{a\in I_1\mid va \in I_1\}$
is also non-zero: for example, $I_1\oo{n} I_1\subseteq I_2$
(it follows from Lemma \ref{lem 5.1.2}(ii) that
 $I_1\oo{\omega } I_1\ne 0$).
One may define
$I_k=\{a\in I_{k-1}\mid va\in I_{k-1}\}$  for every $k> 1$.
By the same way,
$I_m\oo{n} I_k\subseteq I_{k+m}$, so all
$I_k$ are non-zero ideals of $C$.

By Theorem \ref{thm 5.1.1}, there exist an integer $m\ge 0$
and some elements $a_i\in C$, $i=0,\dots, m$, such that
the element
$E_{N,Q} = (1 \otimes e_{NN}) \varphi^{-1}(Q) \in \Cend_{N,Q}$
is presented as
\[
E_{N,Q}  = a_0 + va_1+\dots +v^m a_m,
\]
where $e_{NN}$ is the matrix unit.

Consider an element $0\ne a\in I_{2m+1}$
such that $e_{NN} a e_{NN}  = f(D,v)e_{NN}$,
where $0\ne f(D,v) \in \Cend_1$.
(Note that $e_{NN}I_{2m+1}e_{NN}=0$
is impossible by Lemma \ref{lem 5.1.2}(ii).)
It follows from (\ref{3.1.1}) that
$E_{N,Q}\oo{\omega } I_{2m} \oo{\omega } E_{N,Q}\subseteq C$.
For example,
\[
E_{N,Q}  \oo{0} va \oo{n} E_{N,Q}
 =
 \big( vf_N(v) f(D,v)\oo{n} \varphi^{-1}(f_{N}(v)) \big )
  (1\otimes e_{NN})
\]
lies in $C$
for any $n\in \Zset_+$.
Hence, the conformal algebra $C$ contains matrices of
the form
\[
0\ne x= vg(D,v)(1\otimes e_{NN}),
\]
where $g(D,v)\in \Cend_1$.
By the very same reasons as in the proof of Theorem \ref{thm 5.3.1}
we conclude that
there exists an element $y\in C$
such that
$y\ne 0$
and
$v^ky\in C$ for all $k\ge 0$. Hence,
by Proposition \ref{prop 5.1.3}
we have
$C=\Cend_{N,Q}$.
\end{proof}


\subsection{The case $C\cap vC=0$}\label{subsec5.4}

Let $C\cap vC= 0$ but the sum (\ref{5.1.1}) is still non-direct.
Then there exists a minimal $n\ge 1$ such that
the sum
$C+vC+\dots +v^nC$ is direct.
Let us use the following notation: for a subspace
$X\subseteq \Cend_N$ we denote
$X+vX+\dots +v^k X $
by
$\big(v^{\le k}\big)X$.

The set
$J_1= \{ a\in C
 \mid v^{n+1}a\in \big( v^{\le n} \big )C \}$
 is a non-zero ideal of $C$. If $J_1=C$, then
$\Bbbk[v]C=C\oplus vC \oplus
  \dots \oplus v^nC=\Cend_{N,Q}$.

If $J_1\ne C$, then we may construct
$J_{k+1}=\{a\in J_k\mid v^{n+1}a\in
  \big(v^{\le n}\big)J_k \}$
for $k\ge 1$.
All these $J_k$ are ideals of $C$, and it is clear
that
\begin{equation}\label{5.4.1}
J_l\oo{\omega } J_m \subseteq J_{l+m},
\quad
\big( v^{\le n+k}\big )J_k
 \subseteq \big (v^{\le n+k-j}\big )J_{k-j},\quad
j=1,\dots, k.
\end{equation}
If $J_k\ne 0$ for some $k\ge 1$, then by
Lemma \ref{lem 5.1.2}(ii) we have
$\Bbbk[v]J_k=\Cend_{N,Q}$.
If $J_l \oo{\omega } J_m = 0$,
then $\Bbbk[v]J_l \oo{\omega} J_m
  = \Cend_{N,Q}\oo{\omega } J_m = 0 $
and $J_m=0$.
Since $J_1\ne 0$, the relation (\ref{5.4.1})
 implies $J_k\ne 0 $ for all $k\ge 1$.
It follows from Lemma \ref{lem 5.1.2}(ii)
that
$\Bbbk[v]J_k=\Cend_{N,Q}$ for all $k\ge 1$.
In particular, $\Bbbk[v]J_{n+1}=\Cend_{N,Q}$.
But
(\ref{5.4.1}) and (\ref{3.1.1}) imply that
$C_0=\big(v^{\le n}\big )J_{n+1}$
is a conformal subalgebra of $\Cend_{N,Q}$.
If $n\ge 1$, then $C_0$
  satisfies the conditions of Theorem \ref{thm 5.3.2},
   so
$C_0=\Cend_{N,Q}$.

In any case, we
obtain that there exists a conformal subalgebra
$C\subseteq \Cend_{N,Q}$
such that
$\Cend_{N,Q}=C\oplus vC\oplus \dots \oplus v^{n}C$,
$n\ge 1$.

\begin{thm}\label{thm 5.4.1}
Let $\Cend_{N,Q}=C\oplus vC\oplus \dots \oplus v^{n}C$.
Then $n=0$.
\end{thm}

\begin{proof}
Suppose that $n\ge 1$. For every $a\in C$ we may consider
a unique presentation
\[
v^{(n+1)}a= a_0+va_1 + \dots + v^{(n)} a_n,
 \quad a_i\in C.
\]
Define the map
$\chi:C\to C$ by the rule $\chi(a)=a_n$.
This is an injective
$H$-linear map satisfying the conditions
\begin{equation}\label{5.4.2}
\chi(a\oo{m}x)=\chi(a)\oo{m} x,
\quad
\chi(x\oo{m} a)=x\oo{m}\chi(a) - mx\oo{m-1}a
\end{equation}
for
$a,x\in C$, $m\ge 0$.
Let us consider $\psi:C\to \Cend_{N,Q}$ defined as follows:
\[
\psi(a)=\chi(a) - va, \quad a\in C.
\]
This is an injective $H$-linear map, and
$\psi(C)\cap C=0$. The map $\psi $ satisfies the following
conditions:
\begin{eqnarray}
& \psi(a\oo{m}x)=\psi(a)\oo{m} x,
\quad
\psi(x\oo{m} a)=x\oo{m}\psi(a),
         \label{5.4.3}  \\
& a,x\in C,\quad  m\ge 0.
  \nonumber
\end{eqnarray}
Since
$\Cend_{N,Q}=C\oplus vC\oplus \dots \oplus v^{n}C$,
we may extend $\psi $ to the map
$\bar\psi: \Cend_{N,Q}\to \Cend_{N,Q}$
by the rule
$\bar\psi(v^ka)=v^k \psi(a)$, $k=0,\dots, n$.
It is easy to check that
this $\bar\psi$ satisfies (\ref{5.4.3})
 for every $a\in C$, $x\in \Cend_{N,Q}$:
\begin{eqnarray}
& \bar\psi(a\oo{m}x)=\psi(a)\oo{m} x=a\oo{m} \bar\psi(x),
           \nonumber  \\
& \bar\psi(x\oo{m}a)=\bar\psi(x)\oo{m} a=x\oo{m} \psi(a),
  \nonumber
\end{eqnarray}

In particular, let us put $x=\varphi^{-1}(Q)$,
$\bar \psi(x)=\sum\limits_{t\ge 0} (-D)^{(t)}\otimes X_t$.
Then for an arbitrary element $a\in C$ we have
\[
\bar\psi(x\oo{0}a)= x\oo{0}\psi(a) = \bar\psi(x)\oo{0} a.
\]
Hence,
$Q\psi(a)=X_0a$, and we may conclude that
$\psi(a)=Q^{-1}X_0a$.
Since $\bar \psi $ is defined by
$\bar\psi(v^ka)=v^k\psi(a)$,
$k=0,\dots,n$, $a\in C$,
  we may assume that
$\bar\psi(x)=Q^{-1}X_0x$ for every $x\in \Cend_{N,Q}$.
In particular,
it follows that
$\bar \psi(vx)=v\bar\psi(x)$
for every $x\in \Cend_{N,Q}$.
It is easy to conclude that
\begin{equation}\label{5.4.4}
 \bar\psi(x\oo{m} y) = \bar\psi(x)\oo{m} y = x\oo{m} \bar\psi(y),
 \quad x,y\in \Cend_{N,Q},\ m\ge 0.
\end{equation}
Also, $\bar\psi $ is an injective $H$-linear map.

The only possibility for (\ref{5.4.4}) is
$\bar\psi=\alpha\idd$,
$\alpha\in \Bbbk$. Indeed, let us apply (\ref{5.4.4}) for
$y=\varphi^{-1}(Q)$ and $x=e_{ij}y$.
Denote $B=Q^{-1}X_0\in M_N(\Bbbk(v))$,
and consider
\[
\bar\psi (x)\oo{0}y  =  (1 \otimes Be_{ij}Q)\varphi^{-1}(Q)
 = x\oo{0}\bar\psi(y) = (1\otimes e_{ij}QB)\varphi^{-1}(Q).
\]
In particular,
 we may conclude that
\begin{equation}\label{5.4.5}
Be_{ij}Q=e_{ij}QB,
\end{equation}
so $B$ is a scalar matrix.

Now,
consider
 \begin{eqnarray}
 \bar\psi (x)\oo{1}y
  & = &
   (1\otimes Be_{ij}\partial_v(Q)) \varphi^{-1}(Q)
   =   x\oo{1}\bar\psi(y)
        \nonumber \\
 & = &
 (1 \otimes e_{ij}\partial_v(Q)B ) \varphi^{-1}(Q)
 +
   (1\otimes e_{ij}Q\partial_v(B) ) \varphi^{-1}(Q).
  \label{5.4.6}
\end{eqnarray}
It follows from (\ref{5.4.5}), (\ref{5.4.6}) that
$B e_{ij} \partial_v(Q) = e_{ij}\partial_v(QB) $
for all $i,j=1,\dots, N$.
Since $B$ is a scalar matrix in $M_N(\Bbbk(v))$,
we obtain $B\partial_v(Q) = \partial_v(BQ)$,
 so $B=\alpha I_N$, $\alpha\in \Bbbk$.
Therefore, $\psi(C)\subseteq C$, and
$C\cap vC\ne 0$
in contradiction with $n\ge 1$.
\end{proof}


\subsection{Associative conformal algebras with finite faithful
representation}\label{subsec5.5}

Let us compile the arguments stated above in order to prove the
following

\begin{thm}\label{thm 5.5.1}
Let
 $C\subseteq \Cend_{N}$
 be an irreducible conformal subalgebra. Then either $C$
is a conjugate of $\Curr_N$
or $C=\Cend_{N,Q}$ for
some matrix $Q\in M_N(\Bbbk[v])$, $\det Q\ne 0$.
\end{thm}

\begin{proof}
It follows from Theorem \ref{thm 5.1.1}  that
$\Bbbk[v]C=
\sum\limits_{n\ge 0} v^n C
=\Cend_{N,Q}$ for a suitable matrix~$Q$.
If the sum is direct, then
by Theorem \ref{thm 5.2.6}
$C$ is a conjugate of $\Curr_N$.
If the sum is not direct,
then Theorem \ref{thm 5.4.1} implies $C\cap vC\ne 0$,
so $C=\Cend_{N,Q}$ by Theorem \ref{thm 5.3.2}.
\end{proof}

\begin{lem}\label{lem 5.5.2}
Let $C\ne 0$ be an associative conformal algebra
with a faithful representation of finite type.
If
$\{a\in C \mid C \oo{\omega }a = 0\}
 =\{a\in C\mid a \oo{\omega } C=0\}=0$,
then
$C$
could be embedded into
$\Cend_N$ for some $N$ in such a way that
every proper $C$-submodule of $V_N$ is not faithful.
\end{lem}

\begin{proof}
It is clear that $C\subseteq \Cend_m$ for some $m\ge 1$
(see, e.g., \cite{21}). Consider any descending chain of
faithful $C$-submodules of $V_m$:
\begin{equation}\label{5.5.1}
V_m \supseteq U_1 \supseteq U_2 \supseteq \dots
\end{equation}
It is easy to note (see, e.g., \cite[Lemma~2.1]{1}) that
there exists $n\ge 1$ such that
for every $k\ge 0$ we have
$U_{n+k}\supseteq f_k U_n$ for some $f_k\in H$.

In particular, the $H$-module
$U_n/U_{n+k}$ coincides with its torsion,
so $C\oo{s} U_n \subseteq U_{n+k}$.
Since the $C$-module
$C\oo{\omega } U_n\subseteq V_m$
is faithful, it
is a lower bound of the initial chain (\ref{5.5.1}) in the set of
faithful $C$-submodules of $V_m$. Hence, there exists
a minimal faithful $C$-submodule $U$ of $V_m$.
Since $U$
is finitely generated over $H$ and torsion-free, $U$
is isomorphic to $V_N$ for some $N\ge 1$.
\end{proof}

\begin{thm}\label{thm 5.5.3}
Let $C$ be a conformal algebra with a faithful representation
of finite type. Then:

{\rm (i)} there exists a maximal nilpotent ideal
  ${\mathfrak N}$ of~$C$;

{\rm (ii)} if $C$ is simple, then $C$ is isomorphic to either
$\Curr_N$  or $\Cend_{N,Q}$, $\det Q\ne 0$;

{\rm (iii)} if $C$ is semisimple, then
$C=\bigoplus\limits_{s=1}^{n} I_s$, where $I_s$, $s=1,\dots, n$,
are simple ideals
of~$C$, described in~{\rm (ii)}.
\end{thm}

\begin{proof}
(i) Let $V$ be a faithful $C$-module of finite type.
First, consider the ascending chain of submodules
$U_k = \{u\in V \mid C^k\oo{\omega } u = 0\}$, $k\ge 1$.
Since $V$ is a Noetherian $H$-module, the $H$-module
$V/U_n$ is torsion-free for sufficiently large $n\ge 1$.
In particular, $V/U_n$
is isomorphic as an $H$-module to $V_N$ for a~suitable~$N$.
The ideal $I=\{a\in C\mid C^n\oo{\omega } a  = 0\} \ideal  C$
is nilpotent (or even zero), and $C/I$ acts faithfully on
$V/U_n$.

Therefore, it is sufficient to show that
for every ascending chain
of nilpotent ideals
\begin{equation}\label{5.5.2}
I_1\subseteq I_2 \subseteq \dots ,\quad I_j\ideal C\subseteq \Cend_N,
\ j\ge 1,
\end{equation}
their upper bound $\bigcup\limits_{j\ge 0} I_j$
is also nilpotent.

Let us consider a chain (\ref{5.5.2}) of nilpotent ideals.
Denote by $n_j$ ($j\ge 1$) the nilpotency index of~$I_j$.
Assume that the union of this chain is not nilpotent, so
it could be conjectured that
$n_j<n_{j+1} $ for all $j\ge 1$.

Since the $H$-module $V_N$ is Noetherian, any ascending chain
of submodules becomes stable. In particular,
\begin{eqnarray}
& I_1\oo{\omega }V_N \subseteq I_2\oo{\omega }V_N \subseteq \dots
 = W_1 = I_{j}\oo{\omega } V_N,\quad j\ge l_1 ,
 \nonumber  \\
& I^2_1\oo{\omega }V_N \subseteq I^2_2\oo{\omega }V_N \subseteq \dots
 = W_2 = I^2_{j}\oo{\omega } V_N,\quad j\ge l_2,
 \nonumber  \\
& \hbox to 80mm{\dotfill}   \nonumber  \\
& I^k_1\oo{\omega }V_N \subseteq I^k_2\oo{\omega }V_N \subseteq \dots
 = W_k = I^k_{j}\oo{\omega } V_N,\quad j\ge l_k,
  \nonumber     \\
& \hbox to 80mm{\dotfill}    \nonumber
\end{eqnarray}
It is clear that $W_1\supseteq W_2 \supseteq \dots $,
so there exists an integer $k\ge 1$ such that
$W_k/W_{k+j}$ coincides with its torsion for every
$j\ge 0$ (see, e.g., \cite[Lemma~2.1]{1}).
Hence,
\begin{equation}\label{5.5.3}
C\oo{\omega } W_k \subseteq W_{m}
   \quad \hbox{for all $m>k$}.
\end{equation}
Let us fix some $m>k+1$. Since the sequence $\{n_j\}$
is assumed to be increasing, there exists a number
$p\ge \max\{l_k, l_m\}$
such that $n_p>m$.
Then
$W_k = I_p^k\oo{\omega } V_N$,
$W_m = I_p^m\oo{\omega } V_N$.
The relation (\ref{5.5.3}) implies
\[
I_p^{n_p-m}\oo{\omega } C \oo{\omega } I_p^k\oo{\omega } V_N
   \subseteq I_p^{n_p}\oo{\omega } V_N = 0,
\]
so
\[
0= I_p^{n_p-m}\oo{\omega } C \oo{\omega } I_p^k
  \supseteq I_p^{n_p-m+k+1}.
\]
Therefore, we obtain the contradiction.

(ii) It follows from Lemma \ref{lem 5.5.2} that $C$ has an irreducible
faithful representation of finite type. Theorem \ref{thm 5.5.1} implies
that either $C\simeq \Curr_N$ or $C\simeq \Cend_{N,Q}$, $\det Q\ne 0$.

(iii) By Lemma \ref{lem 5.5.2}, we may assume that $C\subseteq \Cend_N$
and every proper $C$-submodule of $V_N$ is not faithful.

Let $U\subset V_N$ be a maximal $C$-submodule. Then
$I = \Ann_C(U)=\{a\in C \mid a\oo{\omega }U= 0\}\ne 0 $ is an ideal
of~$C$.

Denote by $J=\Ann_C(I)=\{a\in C\mid I\oo{\omega }a = 0\}$
the annihilation ideal of $I$ in $C$. Since
$C$ is semisimple, $I\cap J = 0$.
It is clear that $J\oo{\omega} V_N \subseteq \Ann_{V_N}(I) = U$
since $U$ is a maximal $C$-submodule.

Now, note that
$V_N/U$ is an irreducible $C$-module,
$V_N/U$ is a faithful $I$-module,
$V_N/U$ is a faithful irreducible $C/J$-module.
Since $I$ could be considered as an ideal of $C/J$,
we have $C=I\oplus J$, and $I\simeq C/J$ has a faithful
irreducible representation of finite type.

By (i), $I\simeq \Curr_N$ or $I\simeq \Cend_{N,Q}$, $\det Q\ne 0$.
It is clear that $J$ is also semisimple, and $J$ has
a faithful representation of finite type.

Therefore, every semisimple conformal subalgebra $C$ of
$\Cend_m$  could be presented as
$C=I\oplus J$, where $I$ is simple and $J$ is semisimple,
$I=\Ann_C(J)$. Let us proceed with the decomposition and write
$C=I_1\oplus \dots \oplus I_n \oplus J_n$,
where all $I_j$, $j=1,\dots, n$,
are simple and $J_n$ is semisimple,
$I_1\oplus\dots \oplus I_n =\Ann_C(J_n)$.
Since $\Cend_m $ is a (left and right) Noetherian conformal algebra,
$C$ satisfies a.c.c. for annihilation ideals.
Hence, there exists a sufficiently large integer $n\ge 1$
such that $J_n=0$. Then
$C=I_1\oplus \dots \oplus I_n$.
\end{proof}

\begin{cor}[c.f. \cite{8,21}]\label{cor 5.5.4}
Let $C$ be an associative conformal
algebra of finite type.
If $C$ is simple, then $C\simeq \Curr_N $.
If $C$ is semisimple, then $C$ is a direct sum of simple ideals.
\quad \qed
\end{cor}


\section{Open problems}\label{sec6}

Let us complete the paper with listing some open problems
closely related with the considered one.

\smallskip
{\bf I.}
Describe all irreducible subalgebras of $M_N(W)$
with respect to the canonical action on $(\Bbbk[p])^N$.

Theorems \ref{thm 4.1.2} and \ref{thm 5.1.1} describe the class
of irreducible subalgebras of $M_N(W)$
related with conformal subalgebras of
$\Cend_N$. It is unclear, what would be the answer in
general.

\smallskip
{\bf II.}
Classify irreducible
infinite type
Lie conformal subalgebras of the Lie conformal algebra
$\gc_N=\Cend_N^{(-)}$.

Irreducible Lie conformal subalgebras of finite type
were described in \cite{12} (see
also~\cite{21}).
In \cite{8}, it is conjectured that
the following Lie conformal
subalgebras of $\gc_N$
exhaust all (infinite type)
irreducible Lie conformal subalgebras of $\gc_N$
up to the conjugation by an automorphism of $\Cend_N$
($\det Q\ne 0$ everywhere):
\begin{eqnarray}
& \gc_{N,Q} =\Cend_{N,Q}^{(-)};
                             \nonumber \\
& \oc_{N,Q} = \{a \varphi^{-1}(Q) \mid a\in \Cend_N, \,
 \sigma (a)=-a \},  \quad
   Q^t(-v) = Q(v);
                             \nonumber \\
& \spc_{N,Q} = \{a\varphi^{-1}(Q) \mid a\in \Cend_N, \,
    \sigma (a)= a \},  \quad
     Q^t(-v) = -Q(v).
\nonumber
\end{eqnarray}
Here $\sigma $ is the anti-involution of $\Cend_N$ defined by
the rule  \cite{8}
\[
\sigma(h\otimes A(v))=\sum\limits_{s\ge 0}
 (-D)^{(s)}h\otimes \partial_v^sA^t(v),
\quad h\in H,\  A(v)\in M_N(\Bbbk[v]) .
\]
It is proved in \cite{13} that every irreducible Lie conformal
subalgebra of infinite type
which is an $\Sll_2$-module
(with respect to the action of a Virasoro-like element of $\gc_N$)
is of the type $\oc_{N,Q}$ or $\spc_{N,Q}$.
The result of \cite{29} shows that the conjecture is true for simple
irreducible Lie subalgebras that contain
$\Curr(\Sll_2)$.

\smallskip
{\bf III.}
Describe finitely generated simple and semisimple conformal algebras
of Gel'fand--Kirillov dimension one.

This problem was partially solved in
\cite{23,24,29}. The conjecture stated in the last paper
says that
$\Cend_{N,Q}$, $\det Q\ne 0$,
exhaust all (finitely generated) simple
associative conformal algebras of
Gel'fand--Kirillov dimension one.
We believe it is true, although
 the similar statement does not hold for Lie conformal algebras:
$\gc_{N,Q}$, $\oc_{N,Q}$, $\spc_{N,Q}$,
together with $\Curr {\mathfrak g}$
(where ${\mathfrak g}$ is a simple finitely generated
Lie algebra of Gel'fand--Kirillov dimension one)
is not a complete
list of simple Lie conformal algebras of Gel'fand--Kirillov
dimension one.

\section*{Acknowledgements}
I thank Leonid Bokut and Seok-Jin Kang
for their interest in this work and helpful discussions.
Also, I am grateful to the referee for
many valuable remarks which led to improvement of the presentation.


\end{document}